# RESCALED LOTKA–VOLTERRA MODELS CONVERGE TO SUPER-BROWNIAN MOTION

By J. Theodore Cox[1] and Edwin A. Perkins[2]

*Syracuse University and The University of British Columbia*

We show that a sequence of stochastic spatial Lotka–Volterra models, suitably rescaled in space and time, converges weakly to super-Brownian motion with drift. The result includes both long range and nearest neighbor models, the latter for dimensions three and above. These theorems are special cases of a general convergence theorem for perturbations of the voter model.

**1. Introduction.** In [13], Neuhauser and Pacala introduced a stochastic spatial version of the Lotka–Volterra model for competition between species. We show here that a sequence of these Lotka–Volterra processes, suitably renormalized, converges to super-Brownian motion with a nontrivial drift. We do this by proving a more general convergence theorem, extending the main results of [3] on the voter model. In future work we will show that the above drifts are connected to the questions of co-existence and survival of a rare type in the original Lotka–Volterra model. At present our main results hold for three or more dimensions. Our introduction is structured as follows. In Section 1.1 we describe a special case of the model introduced in [13], and then formulate and state our convergence result. In Section 1.2 we define a class of processes we call *voter model perturbations*, and present a convergence theorem for this class. Our result on Lotka–Volterra models is a special case of this theorem. In Section 1.3 we state and prove a number of corollaries of the main theorem.

1.1. *Lotka–Volterra models.* We suppose that at each site of $\mathbb{Z}^d$ (the $d$-dimensional integer lattice) there is a plant of one of two types. At random times plants die and are replaced by new plants, the times and types

Received October 2003; revised April 2004.
[1]Supported in part by NSF Grant DMS-02-04422.
[2]Supported in part by an NSERC research grant.
*AMS 2000 subject classifications.* Primary 60K35, 60G57; secondary 60F17, 60J80.
*Key words and phrases.* Lotka–Volterra, voter model, super-Brownian motion, spatial competition, coalescing random walk.

depending on the configuration of surrounding plants. The state of the system at time $t$ will be denoted by $\xi_t$, an element of $\{0,1\}^{\mathbb{Z}^d}$, where $\xi_t(x)$ gives the type of the plant at $x$ at time $t$. We have chosen to label the two types 0 and 1; in [13], the types were 1 and 2. To describe the system's evolution, we let $\mathcal{N} \subset \mathbb{Z}^d$ be a finite set not containing the origin, such that $y \in \mathcal{N}$ implies $-y \in \mathcal{N}$. Let $f_i = f_i(\xi) = f_i(x,\xi)$ be the frequency of type $i$ in the neighborhood $x + \mathcal{N}$ in configuration $\xi$,

$$(1.1) \qquad f_i(x,\xi) = \frac{1}{|\mathcal{N}|} \sum_{e \in \mathcal{N}} \mathbb{1}\{\xi(x+e) = i\}, \qquad i = 0, 1.$$

Finally, let $\alpha_0, \alpha_1$ be nonnegative parameters. The dynamics of $\xi_t$ can now be described as follows: at site $x$ in configuration $\xi$, the coordinate $\xi(x)$ makes transitions

$$(1.2) \qquad \begin{aligned} 0 \to 1 &\quad \text{at rate } f_1(f_0 + \alpha_0 f_1), \\ 1 \to 0 &\quad \text{at rate } f_0(f_1 + \alpha_1 f_0). \end{aligned}$$

These rates are interpreted in [13] as follows. A plant of type $i$ dies at rate $f_i + \alpha_i f_{1-i}$, and is replaced by a plant whose type is chosen at random from its neighborhood. In the "death rate" $f_i + \alpha_i f_{1-i}$, $\alpha_i$ measures the strength of inter-specific competition of type $i$, and we have taken the strength of competition due to individuals of the same type to be one. Note that the two configurations, all 0's and all 1's, are both traps. Since $f_0 + f_1 = 1$, the case $\alpha_0 = \alpha_1 = 1$ gives the well-known voter model (see [11] and [3]). In [13], an additional fecundity parameter $\lambda$ allows them to consider populations in which one type has an advantage in replacement. We have chosen to treat only the $\lambda = 1$ case.

Unlike the voter model, the Lotka–Volterra model $\xi_t$ does not have a simple dual process. However, it was shown in [13] that if $\alpha_0 = \alpha_1 = \alpha < 1$, then $\xi_t$ has an *annihilating dual process*, a "double branching annihilating process" in which particles move as random walks, branch, and annihilate each other. Although this process is difficult to analyze, it was instrumental in the proof of Theorem 1 of [13], which states that for $\alpha$ sufficiently small (depending on $\mathcal{N}$, and excluding $\mathcal{N} = \{-1,1\}$ in one dimension), coexistence of types is possible. Here, coexistence means that there is an invariant measure which a.s. concentrates on configurations with infinitely many 0's and infinitely many 1's. On the other hand, comparisons with *biased voter models* (see Section 4) show that for certain values of $(\alpha_0, \alpha_1)$, survival of a given type occurs. More precisely, let $\xi_t^*$ denote the process started from a single 1 at the origin, and 0's everywhere else, and define

$$S = \left\{ (\alpha_0, \alpha_1) : P\left( \sum_{x \in \mathbb{Z}^d} \xi_t^*(x) > 0 \text{ for all } t > 0 \right) > 0 \right\}.$$

Theorem 4 of [13] shows that $\tilde{S} \subset S$, where $\tilde{S}$ is the set of $(\alpha_0, \alpha_1)$ such that

$$(1.3) \qquad 0 \leq \alpha_1 \leq \begin{cases} 1 - \kappa(1 - \alpha_0), & \text{if } 1 - \kappa^{-1} < \alpha_0 < 1, \\ 1 + \kappa^{-1}(\alpha_0 - 1), & \text{if } \alpha_0 > 1, \end{cases}$$

and $\kappa = |\mathcal{N}|$.

We treat here asymptotics for the "low density regime" where there are relatively few plants of one type, which we take to be type 1. It is useful in this context to change our original interpretation, and think now of 0's as representing vacant sites and 1's as representing "particles" which may die or give birth to particles at other sites. We may consider a "measure-valued" version of $\xi_t$ by placing an atom of a given size at each site with a particle. For the voter model case $\alpha_0 = \alpha_1 = 1$, it was shown in [3] (see also [2]) that appropriate low density limits of renormalized voter model processes lead to super-Brownian motion (see Theorem A below). Here we will consider asymptotics for Lotka–Volterra models with the $\alpha_i \to 1$, and will obtain super-Brownian motion with drift in the limit.

Let $\mathcal{M}_f(\mathbb{R}^d)$ denote the space of finite Borel measures on $\mathbb{R}^d$, endowed with the topology of weak convergence of measures. Let $\Omega_{X,D} = D([0,\infty), \mathcal{M}_f(\mathbb{R}^d))$ be the Skorohod space of cadlag $\mathcal{M}_f(\mathbb{R}^d)$-valued paths, and let $\Omega_{X,C}$ be the space of continuous $\mathcal{M}_f(\mathbb{R}^d)$-valued paths with the topology of uniform convergence on compacts. In either case, $X_t$ will denote the coordinate function, $X_t(\omega) = \omega(t)$. Integration of a function $\phi$ with respect to a measure $\mu$ will be denoted by $\mu(\phi)$. For $1 \leq n \leq \infty$, let $C_b^n(\mathbb{R}^d)$ be the space of bounded continuous functions whose partial derivatives of order $n$ or less are also bounded and continuous.

An adapted a.s.-continuous $M_f(\mathbb{R}^d)$-valued process $X_t, t \geq 0$ on a complete filtered probability space $(\Omega, \mathcal{F}, \mathcal{F}_t, P)$ is said to be a *super-Brownian motion with branching rate* $b \geq 0$, *drift* $\theta \in \mathbb{R}$ *and diffusion coefficient* $\sigma^2 > 0$ *starting at* $X_0 \in \mathcal{M}_f(\mathbb{R}^d)$ if it solves the following martingale problem:

(MP) For all $\phi \in C_b^\infty(\mathbb{R}^d)$,

$$(1.4) \qquad M_t(\phi) = X_t(\phi) - X_0(\phi) - \int_0^t X_s\left(\frac{\sigma^2 \Delta \phi}{2}\right) ds - \theta \int_0^t X_s(\phi) \, ds$$

is a continuous $(\mathcal{F}_t)$-martingale, with $M_0(\phi) = 0$ and predictable square function

$$(1.5) \qquad \langle M(\phi) \rangle_t = \int_0^t X_s(b\phi^2) \, ds.$$

The existence and uniqueness in law of a solution to this martingale problem is well known (see, e.g., Theorem II.5.1 and Remark II.5.13 of [14]). Let $P_{X_0}^{b,\theta,\sigma^2}$ denote the law of the solution on $\Omega_{X,C}$ (and also a probability on the space of cadlag paths $\Omega_{X,D}$).

We define our rescaled Lotka–Volterra models following the approach used in [3]. For $N = 1, 2, \ldots$, let $M_N \in \mathbb{N}$ (the set of positive integers), and let $\ell_N = M_N \sqrt{N}$. Let $\mathbf{S_N} = \mathbb{Z}^d/\ell_N$, and let $W_N = (W_N^1, \ldots, W_N^d) \in (\mathbb{Z}^d/M_N \setminus \{0\})$ be a sequence of random vectors such that

(H1)
  (a) $W_N$ and $-W_N$ have the same distribution.
  (b) There is a finite $\sigma^2 > 0$ such that $\lim_{N \to \infty} E(W_N^i W_N^j) = \delta_{ij}\sigma^2$.
  (c) The family $\{|W_N|^2, N \in \mathbb{N}\}$ is uniformly integrable.

Define the kernels $p_N$ by

$$p_N(x) = P\left(\frac{W_N}{\sqrt{N}} = x\right), \qquad x \in \mathbf{S_N}. \tag{1.6}$$

For $\xi \in \{0,1\}^{\mathbf{S_N}}$, define the densities $f_i^N = f_i^N(\xi) = f_i^N(x, \xi)$ by

$$f_i^N(x, \xi) = \sum_{y \in \mathbf{S_N}} p_N(y - x)\mathbb{1}\{\xi(y) = i\}, \qquad i = 0, 1. \tag{1.7}$$

We let $\alpha_i = \alpha_i^N$ depend on $N$, and let $\xi_t^N$ be the process taking values in $\{0,1\}^{\mathbf{S_N}}$ determined by the rates: at site $x$ in configuration $\xi$, the coordinate $\xi(x)$ makes transitions

$$\begin{aligned}
0 \to 1 &\qquad \text{at rate } Nf_1^N(f_0^N + \alpha_0 f_1^N), \\
1 \to 0 &\qquad \text{at rate } Nf_0^N(f_1^N + \alpha_1 f_0^N).
\end{aligned} \tag{1.8}$$

That is, $\xi_t^N$ is the rate-$N$ Lotka–Volterra process determined by the parameters $\alpha_i^N$ (and kernel $p_N$), which we will abbreviate as $LV(\alpha_0^N, \alpha_1^N)$. Note that we recover the original formulation of our process by setting $N = 1$ and letting $W_1$ be uniformly distributed over $\mathcal{N}$, that is, $p_N(x) = \mathbb{1}_{\{x \in \mathcal{N}\}}/|\mathcal{N}|$.

We now consider the measure $X_t^N$ determined by assigning mass $1/N'$ to each site of $\xi_t^N$ with value 1 and mass 0 to all other sites. Here the scaling for the particle mass satisfies $1 \leq N' \leq N$, and will depend on the particular choice of the $W_N$. Given a sequence $N'(N)$, we define the corresponding measure-valued process $X_t^N$ by

$$X_t^N = \frac{1}{N'} \sum_{x \in \mathbf{S_N}} \xi_t^N(x)\delta_x \tag{1.9}$$

($\delta_x$ is the unit point mass at $x$). We make the following assumptions about the initial states $\xi_0^N$:

(H2)
  (a) $\sum_{x \in \mathbf{S_N}} \xi_0^N(x) < \infty$.
  (b) $X_0^N \to X_0 \qquad$ in $\mathcal{M}_f(\mathbb{R}^d)$ as $N \to \infty$.

A consequence of (H2) is that $\sup_N X_0^N(\mathbf{1}) < \infty$, a fact we will frequently use.

*The conditions* (H1) *and* (H2) *will be in force throughout this paper.*

Our basic assumption concerning the rates $\alpha_i^N$ is for $i = 0, 1$,

(H3) $\qquad\qquad \theta_i^N = N(\alpha_i^N - 1) \to \theta_i \in \mathbb{R} \qquad \text{as } N \to \infty.$

We will for the most part focus on Lotka–Volterra models with two types of kernels $p_N$.

(M1) *Long range models.* Let $W_N$ be uniformly distributed on $(\mathbb{Z}^d/M_N) \cap I$, where $I = [-1, 1]^d \setminus \{0\}$, and as $N \to \infty$,

$$\begin{aligned} M_N/\sqrt{N} &\to \infty &&\text{in } d = 1, \\ M_N^2/(\log N) &\to \infty &&\text{in } d = 2, \\ M_N &\to \infty &&\text{in } d \geq 3. \end{aligned}$$

It is simple to check that all the parts of (H1) are satisfied with $\sigma^2 = 1/3$.

(M2) *Fixed kernel models.* Let $M_N \equiv 1$, and let $p(x)$ be an irreducible, symmetric, random walk kernel on $\mathbb{Z}^d$, such that $p(0) = 0$ and $\sum_{x \in \mathbb{Z}^d} x^i x^j p(x) = \delta_{ij} \sigma^2 < \infty$. Define $W_N$ by $P(W_N = x) = p(x)$. It is simple to check that (H1) is satisfied in this case.

As noted before, if we set each $\alpha_i^N = 1$, so that $\theta_0^N = \theta_1^N = 0$, then the $LV(1,1)$ process $\xi_t^N$ is, in fact, the voter model. It was shown in [3] that in this case $X_t^N$ converges weakly in $\Omega_{X,D}$ to super-Brownian motion. More precisely, let $P_N$ denote the law of $X_\cdot^N$. If (M1) holds and $N' = N$, then

(1.10) $\qquad\qquad P_N \Rightarrow P_{X_0}^{2,0,1/3} \qquad \text{as } N \to \infty.$

Under (M2) we have the following (Theorem 1.2 of [3]):

THEOREM A. *Assume* (M2). (a) *If $d \geq 3$ and $N' \equiv N$, then*

$$P_N \Rightarrow P_{X_0}^{2\gamma_e, 0, \sigma^2} \qquad \text{as } N \to \infty.$$

*Here $\gamma_e$ is the "escape probability" of a random walk with step distribution $p$ [see* (1.11) *below].*

(b) *If $d = 2$ and $N' = N/\log N$, then*

$$P_N \Rightarrow P_{X_0}^{4\pi\sigma^2, 0, \sigma^2} \qquad \text{as } N \to \infty.$$

The two-dimensional case in the above theorem is the most delicate and explains why we allowed the possibility of $N' \neq N$ in our definition of $X_t^N$. As explained in [3] (or see Proposition 2.3 below), the voter model may be viewed as a branching random walk with state dependent branching rate $2f_0^N(x, \xi_t^N)$. For $d = 2$, this rate will approach 0 as $N \to \infty$ due to the recurrence of two-dimensional random walk. To counteract this, we increase the branching rate by a factor of $\log N$, or equivalently, reduce the inverse mass per particle by a factor of $\log N$. As we will only treat either the fixed kernel case with $d \geq 3$ or the long range case below, we will assume that

$$N' = N \quad \text{in the rest of this work.}$$

Let us return now to the Lotka–Volterra models $\xi_t^N$. We let $P_N$ denote the law of $X_\cdot^N = \frac{1}{N} \sum_{x \in \mathbf{S_N}} \xi_\cdot^N(x) \delta_x$ on $\Omega_{X,D}$. Under the assumption (H3) on the rates $\alpha_i^N$, we again have convergence to super-Brownian motion, but this time with a (possibly) nonzero drift. Recall that (H1) and (H2) are always in force.

THEOREM 1.1. *Assume* (H3) *and* (M1). *Then* $P_N \Rightarrow P_{X_0}^{2, -\theta_1, 1/3}$ *as* $N \to \infty$.

Next, we consider the fixed kernel case (M2). This time, to specify the parameters in the limiting super-Brownian motion, we must introduce a coalescing random walk system $\{\hat{B}_t^x, x \in \mathbb{Z}^d\}$. Each $\hat{B}_t^x$ is a rate 1 random walk on $\mathbb{Z}^d$ with kernel $p$, with $\hat{B}_0^x = x$. The walks move independently until they collide, and then move together after that. For finite $A \subset \mathbb{Z}^d$, let $\tau(A) = \inf\{s : |\{\hat{B}_s^x, x \in A\}| = 1\}$ be the time at which the particles starting from $A$ coalesce into a single particle, and write $\tau(a, b, \dots)$ when $A = \{a, b, \dots\}$. For $d \geq 3$, define the "escape" probability (used in Theorem A) by

$$\gamma_e = \sum_{e \in \mathbb{Z}^d} p(e) P(\tau(0, e) = \infty). \tag{1.11}$$

Note that $\gamma_e$ is the probability that a discrete time random walk with step distribution $p$, starting at the origin, never returns to the origin. We also define

$$\begin{aligned} \beta &= \sum_{e, e' \in \mathbb{Z}^d} p(e) p(e') P(\tau(e, e') < \infty, \tau(0, e) = \tau(0, e') = \infty), \\ \delta &= \sum_{e, e' \in \mathbb{Z}^d} p(e) p(e') P(\tau(0, e) = \tau(0, e') = \infty). \end{aligned} \tag{1.12}$$

Here we are considering a system of 3 coalescing random walks starting at $0$, $e$ and $e'$, where $e$ and $e'$ are independent with law $p$. Then $\beta$ is the probability the walks starting at $e$ and $e'$ coalesce, but this coalescing system does not meet the random walk starting at 0, while $\delta$ is the strictly larger probability that the coalescing system starting at $\{e, e'\}$ does not meet the random walk starting at 0.

THEOREM 1.2. *Assume* (H3), (M2) *and* $d \geq 3$. *Then* $P_N \Rightarrow P_{X_0}^{2\gamma_e, \theta, \sigma^2}$ *as* $N \to \infty$, *where* $\theta = \theta_0 \beta - \theta_1 \delta$.

Although Theorem 1.1 is a simpler result than Theorem 1.2, it includes the low-dimensional case $d \leq 2$. Theorem A suggests that it should be possible to extend Theorem 1.2 to the more delicate two-dimensional setting, with $N' = N/\log N$ and a different drift arising from asymptotic versions of $\beta$ and $\delta$. This is the objective of parallel work.

In Theorem 1.1 there is no $\theta_0$ dependence in the limiting law. This suggests the possibility of a long range limit theorem without insisting that $\alpha_0^N$ approach 1. This is, indeed, the case and in a forthcoming paper we will establish a long range limit theorem for fixed $\alpha_0 \in [0,1]$ and $\alpha_1^N$ as above. The argument, based on a combination of ideas used here and in the corresponding convergence for the long range contact process [6], suggests that a unification and generalization of these results should be possible.

Our motivation for this work is two-fold. First, it has been shown in recent years that a number of different spatial stochastic systems at or near criticality, and above a "critical dimension," converge to super-Brownian motion or a near relative when suitably rescaled. This includes lattice trees above 8 dimensions [4], long-range contact processes above 1 dimension [6], oriented percolation above 4 spatial dimensions [9] and, of course, the voter model (Theorem A above). (See [15] for a nice survey.) It is natural to ask if the same is true for the $LV(\alpha_0, \alpha_1)$ models. The above results are steps in this direction, but, more generally, one could ask if such a limit theorem will hold [in the context of (M2)] with zero limiting drift for any "critical" $LV(\alpha_0, \alpha_1)$ model. (Of course, one must define "critical" here.) A second motivation for proving any limit theorem is to actually use it to study the more complicated approximating systems—especially, as is the case here, when there are few tools available for their study. In a forthcoming paper we will use Theorem 1.2 to refine the survival and co-existence results of [13] mentioned earlier for $(\alpha_0, \alpha_1)$ near $(1, 1)$.

1.2. *Voter model perturbations.* In view of assumption (H3), the Lotka–Volterra models $\xi_t^N$ can be viewed as small *perturbations* of the voter model. To see this, we first rewrite the rates in (1.8) in the form

$$
\begin{aligned}
0 \to 1 & \quad \text{at rate } N f_1^N + \theta_0^N (f_1^N)^2, \\
1 \to 0 & \quad \text{at rate } N f_0^N + \theta_1^N (f_0^N)^2.
\end{aligned}
\tag{1.13}
$$

Adopting the notation of [11], the Lotka–Volterra model $\xi_t^N$ is the *spin-flip system* with rate function $c_N(x, \xi)$ [which gives the rate at which coordinate $\xi(x)$ changes to $1 - \xi(x)$],

$$c_N(x, \xi) = N c_N^v(x, \xi) + c_N^*(x, \xi), \tag{1.14}$$

where $c_N^v(x,\xi)$ is the voter model rate function

$$(1.15) \qquad c_N^v(x,\xi) = \sum_{e \in \mathbf{S_N}} p_N(e)\mathbb{1}\{\xi(x+e) \neq \xi(x)\}$$

and $c_N^*(x,\xi)$ is the "perturbation"

$$(1.16) \quad c_N^*(x,\xi) = \theta_0^N(f_1^N(x,\xi))^2\mathbb{1}\{\xi(x)=0\} + \theta_1^N(f_0^N(x,\xi))^2\mathbb{1}\{\xi(x)=1\}.$$

We will generalize the above, defining a wider class of voter model perturbations, and prove convergence to super-Brownian motion for these processes (hence, including Theorems 1.1 and 1.2 as special cases). First, we need some additional notation. Let $P_F$ denote the set of finite subsets of $\mathbb{Z}^d$. For $A \in P_F$, $x \in \mathbf{S_N}$, $\xi \in \{0,1\}^{\mathbf{S_N}}$, define

$$\chi_N(A,x,\xi) = \prod_{e \in A/\ell_N} \xi(x+e).$$

We assume now that $c_N(x,\xi)$ is a function of the form given in (1.14), where $c_N^v(x,\xi)$ is as in (1.15), and $c_N^*(x,\xi)$ is given by

$$(1.17) \quad c_N^*(x,\xi) = \sum_{A \in P_F} \chi_N(A,x,\xi)(\beta_N(A)\mathbb{1}\{\xi(x)=0\} + \delta_N(A)\mathbb{1}\{\xi(x)=1\}).$$

Here $\beta_N$ and $\delta_N$ are real-valued functions on $P_F$ (which may take negative values), but we will assume throughout that

$$(1.18) \qquad c_N(x,\xi) \geq 0 \qquad \text{for all } x, \xi.$$

It is easy to check that the Lotka–Volterra rates can be written as in (1.17) [see (1.25) and (1.26) below].

We now make a number of assumptions on the kernels $p_N$ and on the perturbation rates $\beta_N$ and $\delta_N$.

*Kernel assumptions.* The kernel assumptions (K1)–(K3) below are similar to the ones in [3]. We assume that the $p_N$ are given by (1.6) [recall (H1) is in force], and we let $\{\hat{B}_t^{N,x}, x \in \mathbf{S_N}\}$ denote a rate-$N$ continuous time coalescing random walk system on $\mathbf{S_N}$ with step distribution $p_N$ such that $\hat{B}_0^{N,x} = x$. For finite $A \subset \mathbf{S_N}$, let $\hat{\tau}^N(A)$ denote the time at which all particles starting from $A$ have coalesced into a single particle,

$$\hat{\tau}^N(A) = \inf\{t \geq 0 : |\{\hat{B}_t^{N,x}, x \in A\}| = 1\}.$$

We will also need a collection of independent (noncoalescing) rate-$N$ continuous time random walks with step distribution $p_N$, which we will denote $\{B_t^{N,x} : x \in \mathbf{S_N}\}$, such that $B_0^{N,x} = x$. We can now state the kernel assumptions.

We assume there is a constant $\gamma \geq 0$ and a positive sequence $\{\varepsilon_N^*\}$ with $\varepsilon_N^* \to 0$ and $N\varepsilon_N^* \to \infty$ as $N \to \infty$, such that the following hold:

(K1) $$\lim_{N \to \infty} NP(B_{\varepsilon_N^*}^{N,0} = 0) = 0.$$

(K2)
$$\lim_{N \to \infty} \sum_{e \in \mathbf{S_N}} p_N(e) P(\hat{\tau}^N(\{0,e\}) \in (\varepsilon_N^*, t]) = 0 \qquad \text{for all } t > 0,$$
$$\lim_{N \to \infty} \sum_{e \in \mathbf{S_N}} p_N(e) P(\hat{\tau}^N(\{0,e\}) > \varepsilon_N^*) = \gamma.$$

For $A \in P_F$, let $\tau^N(A) = \hat{\tau}^N(A/\ell_N)$, and put $\sigma_N(A) = P(\tau^N(A) \leq \varepsilon_N^*)$. [We make the convention $\tau^N(\varnothing) = 0$, so $\sigma_N(\varnothing) = 1$.] The last kernel assumption we need is

(K3) $$\sigma(A) = \lim_{N \to \infty} \sigma_N(A) \qquad \text{exists for all } A \in P_F.$$

We ask the reader to distinguish between the function $\sigma(\cdot)$ defined above and the variance parameter $\sigma^2$ in (H1).

We will see below that the conditions (K1)–(K3) hold if the kernels $p_N$ are either of the long range (M1) or fixed kernel (M2) type.

A key step will be to show that local spatial averages of microscopic quantities like the local density of 1's or 0's near a 1 converge to certain coalescing probabilities (like $\beta$ or $\delta$) as $N \to \infty$. The spatial averaging will be implemented by taking a conditional expectation with respect to the process up to time $t - \varepsilon_N^*$, where $t$ is the current time. So $\varepsilon_N^*$ must be large enough to allow enough time for the averaging [hence, (K1) and (K2)], but still approach 0 to ensure locality of the averaging.

*Perturbation assumptions.* We may assume without loss of generality that
$$\beta_N(A) = \delta_N(A) = 0 \qquad \text{if } 0 \in A.$$

To see why this is the case, note that the value of $\beta_N(A)$ is irrelevant when $0 \in A$ because $\chi_N(A, x, \eta)\mathbb{1}(\eta(x) = 0) = 0$. If we define
$$\delta_N'(A) = \begin{cases} 0, & \text{if } 0 \in A, \\ \delta_N(A) + \delta_N(A \cup \{0\}), & \text{if } 0 \notin A, \end{cases}$$

then a short calculation shows that replacing $\delta_N$ with $\delta_N'$ does not change $c_N^*(x, \eta)$.

The assumptions we now make appear somewhat technical, but in Section 1.3 we will show that they can be simplified (or hold automatically) in some natural special cases. Roughly speaking, (P1) says that the "perturbations" $\beta_N$ and $\delta_N$ are appropriately bounded, (P2) and (P3) say that these rates converge in a well-behaved way, and we require (P4) and (P5) in order to make comparisons with the biased voter model in Section 4. As usual, $\ell_1(P_F)$ is the space of functions $f : P_F \to \mathbb{R}$ such that $\|f\|_1 = \sum_{A \in P_F} |f(A)| < \infty$.

(P1) $$\sup_N \sum_{A \in P_F} \max(|A|, 1)(|\beta_N(A)| + |\delta_N(A)|) < \infty.$$

(P2) There exist functions $\beta, \delta$ on $P_F$ such that
$$\beta_N \to \beta \quad \text{and} \quad \delta_N \to \delta \quad \text{pointwise on } P_F \text{ as } N \to \infty.$$

(P3) If $\sigma(\cdot)$ is in (K3), then as $N \to \infty$,
$$\beta_N(\cdot)\sigma_N(\cdot) \to \beta(\cdot)\sigma(\cdot) \quad \text{and} \quad \delta_N(\cdot)\sigma_N(\cdot \cup \{0\}) \to \delta(\cdot)\sigma(\cdot \cup \{0\})$$
in $\ell_1(P_F)$.

(P4) There is a constant $k_\delta > 0$ such that for all $\xi \in \{0, 1\}^{\mathbb{Z}^d}$ with $\xi(0) = 1$,
$$\sum_{A \in P_F} \delta_N(A) \prod_{a \in A} \xi(a) \geq -k_\delta \sum_{y \in \mathbb{Z}^d} p_N(y/\ell_N)(1 - \xi(y)).$$

(P5) $\beta_N(\varnothing) = 0$.

Condition (P1) and (1.18) imply that the rates $c_N(x, \eta)$ above determine a unique $\{0, 1\}^{\mathbf{S_N}}$-valued Feller process. More specifically, consider the associated Markov pregenerator

(1.19) $$\Omega_N f(\xi) = \sum_{x \in \mathbf{S_N}} c_N(x, \xi)(f(\xi^x) - f(\xi)),$$

defined for functions $f : \mathbf{S_N} \to \mathbb{R}$ which depend on only finitely many coordinates. Here $\xi^x$ is the configuration $\xi$ with the coordinate at $x$ flipped to $1 - \xi(x)$. It is straightforward to check that (P1) and (1.18) imply the hypotheses of Theorem B3 of [12], and so there is a unique Feller process $\xi_\cdot^N$ whose generator is the closure of $\Omega_N$.

For our main result, Theorem 1.3, we assume now that the conditions (1.18), (H1), (H2), (K1)–(K3) and (P1)–(P5) hold, and $\xi_\cdot^N$ is the corresponding voter model perturbation. As before, $X_\cdot^N$ is the measure-valued process determined by $\xi_\cdot^N$, $X_t^N = (1/N) \sum_{x \in \mathbf{S_N}} \xi_t^N(x) \delta_x$, and $P_N$ is the law of $X_\cdot^N$ on $\Omega_{X,D}$.

THEOREM 1.3. *As $N \to \infty$, $P_N \Rightarrow P_{X_0}^{2\gamma, \theta, \sigma^2}$, where $\gamma$ is given in* (K2),

(1.20) $$\theta = \sum_{A \in P_F} \beta(A)\sigma(A) - \sum_{A \in P_F} (\beta(A) + \delta(A))\sigma(A \cup \{0\}),$$

*and $\sigma(\cdot)$ is given in* (K3).

REMARK 1.4. Our assumption that $\beta_N(A) = \delta_N(A) = 0$ if $0 \in A$ implies that $\beta(A) = \delta(A) = 0$ if $0 \in A$. Therefore, letting $P_F' = \{A \in P_F : 0 \notin A\}$, the sums over $P_F$ in (1.20) can be replaced by sums over $P_F'$. Similarly, in (P3), we need only consider convergence in $\ell_1(P_F')$.

1.3. *Applications of Theorem* 1.3. In this section we specialize Theorem 1.3 to kernels $p_N$ which satisfy (M1) or (M2). We will see that in each case, the kernel conditions (K1)–(K3) hold, and that some of the perturbation conditions may be simplified. We also show that the Lotka–Volterra Theorems 1.1 and 1.2 follow from Theorem 1.3. We consider first the fixed kernel case.

Assume first that (M2) holds [and, hence, (H1)], and $d \geq 3$. Then the conditions (K1)–(K3) follow for any sequence $\varepsilon_N^* \to 0$ such that $\varepsilon_N^* \gg N^{-1/3}$. To check (K1), we make use of the local limit theorem bound (see Lemma A.3 of [3], e.g.), $P(B_t^0 = 0) \leq C t^{-d/2}$ for some constant $C$. Since $d \geq 3$,

$$NP(B_{\varepsilon_N^*}^{N,0} = 0) = NP(B_{N\varepsilon_N^*}^0 = 0) \leq C(N\varepsilon_N^{*3})^{-1/2} \to 0 \qquad \text{as } N \to \infty.$$

Next,

$$\sum_{e \in \mathbf{S_N}} p_N(e) P(\hat{\tau}_N(0, e) > \varepsilon_N^*) = \sum_{e \in \mathbb{Z}^d} p(e) P(\tau(0, e) > N\varepsilon_N^*)$$

$$\to \sum_{e \in \mathbb{Z}^d} p(e) P(\tau(0, e) = \infty) = \gamma_e.$$

A similar calculation, using transience of the random walks, shows that the first limit in (K2) holds. For $A \in P_F$,

$$\sigma_N(A) = P(\tau^N(A) \leq \varepsilon_N^*) = P(\tau(A) \leq N\varepsilon_N^*) \to P(\tau(A) < \infty) = \sigma(A),$$

so (K3) holds as well. Furthermore, a little rearrangement shows that we may rewrite the limiting drift $\theta$ given in (1.20) in Theorem 1.3 in the form

$$\theta = \sum_{A \in P_F} \beta(A) P(\tau(A) < \infty, \tau(A \cup \{0\}) = \infty)$$

(1.21)

$$- \sum_{A \in P_F} \delta(A) P(\tau(A \cup \{0\}) < \infty).$$

We can now present several corollaries of Theorem 1.3. We will assume, of course, that the rates $c_N(x, \xi)$ are nonnegative and are given by (1.14) and (1.17), and that (H2) and (M2) hold, and $d \geq 3$, but all other assumptions will be specified. We will consider the alternative conditions

(P1)′ $\qquad \beta_N(A) = \delta_N(A) = 0 \qquad$ if $|A| > n_0$ for some finite $n_0$,

and for some $\beta, \delta \in \ell_1(P_F)$,

(P3)′ $\qquad\qquad\qquad \beta_N \to \beta \quad \text{and} \quad \delta_N \to \delta \qquad \text{in } \ell_1(P_F).$

COROLLARY 1.5. *Assume that the perturbation rates $\{\beta_N\}$, $\{\delta_N\}$ satisfy (P1), (P3)$'$, (P4) and (P5). Then $P_N \Rightarrow P_{X_0}^{2\gamma_e,\theta,\sigma^2}$ as $N \to \infty$, where $\gamma_e$ is the escape probability in (1.11) and $\theta$ is the drift specified in (1.21).*

PROOF. To apply Theorem 1.3, it suffices to check that (P2) and (P3) hold. It is clear that (P3)$'$ implies (P2), and an easy uniform integrability argument using $\sigma_N \le 1$ shows that (P3)$'$ also implies (P3) [recall (K3)]. Thus, the conclusion of Theorem 1.3 holds. □

COROLLARY 1.6. *Assume that the perturbation rates $\{\beta_N\}$, $\{\delta_N\}$ satisfy (P1)$'$, (P3)$'$, (P4) and (P5). Then $P_N \Rightarrow P_{X_0}^{2\gamma_e,\theta,\sigma^2}$ as $N \to \infty$, where $\gamma_e$ is the escape probability in (1.11), and $\theta$ is the drift specified in (1.21).*

PROOF. It is easy to check that (P1)$'$ and (P3)$'$ imply (P1), so we may apply Corollary 1.5. □

If we consider kernels $p$ with finite range (as for simple symmetric random walk), then the technical condition (P4) follows automatically from (a weaker version of) (P1).

LEMMA 1.7. *Assume (M2) and that $p$ has finite range. If*

$$(1.22) \qquad \sup_N \sum_{A \in P_F} \delta_N(A)^- < \infty,$$

*then (P4) holds.*

PROOF. The fact that $c_N(x,\xi) \ge 0$ implies that if $\xi \in \{0,1\}^{\mathbb{Z}^d}$ and $\xi(0) = 1$, then

$$\sum_{A \subset \mathbb{Z}^d} \delta_N(A) \prod_{a \in A} \xi(a) \ge -N \sum_{y \in \mathbb{Z}^d} p(y)(1 - \xi(y)) = -N f_0(0,\xi),$$

where $f_0(x,\xi) = \sum_y p(y-x)(1-\xi(y))$. If $f_0(0,\xi) = 0$, then (P4) holds trivially by the above. If $f_0(0,\xi) > 0$, then the finite range assumption implies that for some $\varepsilon > 0$, $f_0(0,\xi) \ge \varepsilon$. Then (1.22) implies that for some $C > 0$,

$$\sum_{A \in P_F} \delta_N(A) \prod_{a \in A} \xi(a) \ge - \sum_{A \in P_F} \delta_N(A)^- \ge -C.$$

Since $f_0(0,\xi) \ge \varepsilon$, $-C \ge -(C/\varepsilon) f_0(\xi)$, and (P4) follows in this case as well. □

COROLLARY 1.8. *Assume that the perturbation rates $\{\beta_N\}$, $\{\delta_N\}$ satisfy (P1)$'$, (P3)$'$ and (P5), and $p$ has finite range. Then $P_N \Rightarrow P_{X_0}^{2\gamma_e,\theta,\sigma^2}$ as $N \to \infty$, where $\gamma_e$ is the escape probability in (1.11), and $\theta$ is given in (1.21).*

PROOF. By Lemma 1.7, (P4) holds, and so the result is immediate from the previous corollary. □

We consider now the long range case, and will suppose that (M1) [and, hence, (H1)] hold until further notice. To verify that the kernel conditions (K1)–(K3) hold for suitable $\varepsilon_N^*$ and $\sigma(A)$, we rely on results from [3].

The first fact we need is that

$$\lim_{N\to\infty} \sup_{A\in P_F, |A|\geq 2} P(\tau^N(A) \leq t) = 0 \qquad \text{for all } t \geq 0. \tag{1.23}$$

To prove this, we need only take the sup over $|A| = 2$ in the above, but this case is covered in the proof of Theorem 5.1(a) of [3]. Only minor notational changes in that argument are required. We also need Lemma 5.2 of [3], which states that there is a finite constant $C$ such that for all $t \geq 0$,

$$P(B_t^{N,0} = 0) \leq \exp\left(\frac{-Nt}{2}\right) + \frac{C}{M_N^d(Nt+1)^{d/2}}.$$

The condition (K1) follows easily from this last estimate for any $\varepsilon_N^* \to 0$, provided $\varepsilon_N^* \gg N^{-1/3}$ for $d \geq 3$, $\varepsilon_N^* \gg \max(M_N^{-2}, 4\log N/N)$ for $d = 2$, and $\varepsilon_N^* \gg \max(NM_N^{-2}, 4\log N/N)$ for $d = 1$. If we set $\gamma = 1$, then the kernel condition (K2), for any sequence $\varepsilon_N^* \to 0$, is an immediate consequence of (1.23). Setting $\sigma(A) = \mathbb{1}\{|A| \leq 1\}$, condition (K3) also follows from (1.23). In view of the above Remark 1.4, the drift $\theta$ in Theorem 1.3 takes the form

$$\theta = \left[\sum_{a\in\mathbb{Z}^d} \beta(\{a\})\right] - \delta(\varnothing). \tag{1.24}$$

As in the fixed kernel case, we consider two alternative perturbation assumptions:

(P1)″ $$\sup_N \sum_A (|\beta_N(A)| + |\delta_N(A)|) < \infty,$$

(P3)″ $$\{\beta_N(\{a\})\}_{a\in\mathbb{Z}^d} \to \{\beta(\{a\})\}_{a\in\mathbb{Z}^d} \qquad \text{in } \ell_1(\mathbb{Z}^d).$$

Recall that we are assuming (H2) and (M1).

COROLLARY 1.9. *Assume that the perturbation rates $\{\beta_N\}$, $\{\delta_N\}$ satisfy* (P1)′, (P1)″, (P2), (P3)″, (P4) *and* (P5). *Then $P_N \Rightarrow P_{X_0}^{2,\theta,1/3}$ as $N \to \infty$, where $\theta$ is given in* (1.24).

PROOF. To apply Theorem 1.3, we need only check that (P1) and (P3) hold. Condition (P1) is immediate from (P1)′ and (P1)″. For (P3), we note by (1.23) that there is a sequence $\eta_N \to 0$ as $N \to \infty$ such that

$$\sum_{A \in P_F, A \neq \varnothing} |\delta_N(A)| \sigma_N(A \cup \{0\}) = \sum_{A \in P'_F, A \neq \varnothing} |\delta_N(A)| \sigma_N(A \cup \{0\})$$

$$\leq \eta_N \sum_{A \in P'_F, A \neq \varnothing} |\delta_N(A)|$$

$$\leq \eta_N C \to 0,$$

the last inequality by (P1)″. A similar argument shows that

$$\lim_{N \to \infty} \sum_{A \in P_F, |A| > 1} |\beta_N(A)| \sigma_N(A) = 0.$$

These last two results, (P3)″ and $\lim_{N \to \infty} \delta_N(\varnothing) = \delta(\varnothing)$ [which follows from (P2)] imply (P3), so we are done. □

We now derive Theorems 1.1 and 1.2 as applications of Corollary 1.9 and Corollary 1.6, respectively.

PROOF OF THEOREMS 1.1 AND 1.2. As previously noted, the rate function $c_N(x,\xi)$ for the Lotka–Volterra rates (1.8) can be written in the form $Nc_N^v(x,\xi) + c_N^*(x,\xi)$, where $c_N^v(x,\xi)$ is given in (1.15) and $c_N^*(x,\xi)$ is given in (1.16). For configurations $\xi$ with $\xi(x) = 1$, one can rewrite (1.16) in the form

$$\theta_1^N - 2\theta_1^N \sum_{e \in \mathbf{S_N}} p_N(e) \xi(x+e) + \theta_1^N \sum_{e,e' \in \mathbf{S_N}} p_N(e) p_N(e') \xi(x+e) \xi(x+e').$$

It follows easily that if we define $\beta_N$ and $\delta_N$ by

$$(1.25) \quad \beta_N(A) = \begin{cases} \theta_0^N (p_N(a/\ell_N))^2, & A = \{a\}, \\ 2\theta_0^N p_N(a/\ell_N) p_N(a'/\ell_N), & A = \{a, a'\}, \\ 0, & \text{otherwise,} \end{cases}$$

and

$$(1.26) \quad \delta_N(A) = \begin{cases} \theta_1^N, & A = \varnothing, \\ \theta_1^N [(p_N(a/\ell_N))^2 - 2p_N(a/\ell_N)], & A = \{a\}, \\ 2\theta_1^N p_N(a/\ell_N) p_N(a'/\ell_N), & A = \{a, a'\}, \\ 0, & \text{otherwise,} \end{cases}$$

then (1.17) is satisfied.

Before considering the two types of models separately, we note that condition (P4) is satisfied in both cases. This is because (1.13) shows that for $\xi \in \{0,1\}^{\mathbf{S_N}}$ with $\xi(x) = 1$,

$$\sum_{A \subset \mathbb{Z}^d} \delta_N(A) \chi_N(A, x, \xi) = \theta_1^N (f_0^N(x, \xi))^2 \geq -|\theta_1^N| f_0^N(x, \xi).$$

This implies that for $\xi \in \{0,1\}^{\mathbb{Z}^d}$ with $\xi(0) = 1$,

$$\sum_{A \subset \mathbb{Z}^d} \delta_N(A) \prod_{a \in A} \xi(a) \geq -|\theta_1^N| \sum_{y \in \mathbb{Z}^d} p_N(y/\ell_N) \mathbb{1}\{\xi(y) = 0\},$$

and, thus, (P4) follows. Conditions (P1)$'$ (with $n_0 = 2$) and (P5) are also clear for both models.

Consider the long range model (M1), and let $\Gamma_N = ([-M_N, M_N]^d \cap \mathbb{Z}^d) \setminus \{0\}$. The formulas for $\beta_N$ and $\delta_N$ simplify to

$$\beta_N(A) = \begin{cases} \theta_0^N \mathbb{1}\{a \in \Gamma_N\}/|\Gamma_N|^2, & A = \{a\}, \\ 2\theta_0^N \mathbb{1}\{a, a' \in \Gamma_N\}/|\Gamma_N|^2, & A = \{a, a'\}, \\ 0, & |A| \neq 1 \text{ or } 2, \end{cases}$$

and

$$\delta_N(A) = \begin{cases} \theta_1^N, & A = \varnothing, \\ \theta_1^N \mathbb{1}\{a \in \Gamma_N\} \left[\frac{1}{|\Gamma_N|^2} - \frac{2}{|\Gamma_N|}\right], & A = \{a\}, \\ 2\theta_1^N \mathbb{1}\{a, a' \in \Gamma_N\}/|\Gamma_N|^2, & A = \{a, a'\}, \\ 0, & |A| > 2. \end{cases}$$

If we set $\beta(A) = 0$ for all $A$, $\delta(\varnothing) = \theta_1$ and $\delta(A) = 0$ for $A \neq \varnothing$, then clearly (P2) holds. It is also trivial now to verify (P1)$''$ and (P3)$''$. Theorem 1.1 is thus a consequence of Corollary 1.9.

Consider now the fixed kernel model (M2). Due to the assumption $p_N(a/\ell_N) = p(a)$, $\beta_N$ and $\delta_N$ only depend on $N$ through $\theta_i^N$. Therefore, if we define $\beta(A)$ and $\delta(A)$ as $\beta_N(A)$ and $\delta_N(A)$, but with $\theta_i$ in place of $\theta_i^N$, (P3)$'$ is a simple consequence of (H3). The hypotheses of Corollary 1.6 are therefore valid.

It remains only to verify the form of the drift $\theta$ given in Corollary 1.6. Recall the definitions of $\beta$ and $\delta$ from (1.12). The term involving the $\beta(A)$'s in the drift $\theta$ of (1.21) equals

$$\sum_A \beta(A) P(\tau(A) < \infty, \tau(A \cup \{0\}) = \infty)$$

$$= \theta_0 \sum_e p^2(e) P(\tau(0, e) = \infty)$$

$$+ \theta_0 \sum_{e \neq e'} p(e) p(e') P(\tau(e, e') < \infty, \tau(0, e, e') = \infty)$$

$$= \theta_0 \sum_{e, e'} p(e) p(e') P(\tau(e, e') < \infty, \tau(0, e) = \tau(0, e') = \infty) = \theta_0 \beta.$$

The term involving the $\delta(A)$'s is

$$\theta_1 \left[1 + \sum_e (p(e)^2 - 2p(e)) P(\tau(0, e) < \infty) + \sum_{e \neq e'} p(e) p(e') P(\tau(0, e, e') < \infty)\right]$$

$$= \theta_1 \Bigg[ 1 + \sum_{e,e'} p(e)p(e')(1 - P(\tau(0,e,e') = \infty))$$

$$- 2\sum_e p(e)P(\tau(0,e) < \infty) \Bigg]$$

$$= \theta_1 \Bigg[ 2\sum_e p(e)P(\tau(0,e) = \infty)$$

$$- \sum_{e,e'} p(e)p(e')(P(\tau(0,e) = \infty) + P(\tau(0,e) < \infty, \tau(0,e') = \infty)) \Bigg]$$

$$= \theta_1 \Bigg[ \sum_e p(e)P(\tau(0,e) = \infty)$$

$$- \sum_{e,e'} p(e)p(e')P(\tau(0,e') < \infty, \tau(0,e) = \infty) \Bigg] = \theta_1 \delta.$$

In the next to last line we used symmetry to interchange $e$ and $e'$. This shows the drift in Corollary 1.6 equals that in Theorem 1.2, and so Theorem 1.2 is proved as well. □

For our final application of Theorem 1.3, we consider rescaled Lotka–Volterra models in which the *dispersion* kernel is still $p_N$, but the *competition* kernels for the two types may be different. We focus on the fixed kernel case (M2) with $d \geq 3$, and fix a pair of competition kernels $p^b$ and $p^d$ on $\mathbb{Z}^d$. The latter two kernels are arbitrary laws on $\mathbb{Z}^d$ satisfying $p^b(0) = p^d(0) = 0$, while the dispersal kernel $p$ still is as in (M2). The rates for the rescaled process $\xi_t^N$ on $\mathbf{S_N} = \mathbb{Z}^d / \sqrt{N}$ are now given by

$$(1.27) \qquad \begin{aligned} 0 \to 1 & \quad \text{at rate } N f_1^N (f_0^{b,N} + \alpha_0^N f_1^{b,N}), \\ 1 \to 0 & \quad \text{at rate } N f_0^N (f_1^{d,N} + \alpha_1^N f_0^{d,N}). \end{aligned}$$

Here $f_i^{b,N}$ is the local density of type $i$ with respect to the rescaled kernel $p_N^b$, and similarly for $f_i^{d,N}$. We continue to assume (H2) and (H3). As before, $X_t^N$ is the empirical measure which assigns mass $1/N$ to the site of each 1 in $\xi_t^N$, and $P_N$ is its law. Finally, we define

$$\beta' = \sum_{e,e' \in \mathbb{Z}^d} p(e)p^b(e')P(\tau(e,e') < \infty, \tau(0,e) = \tau(0,e') = \infty),$$

$$\delta' = \sum_{e,e' \in \mathbb{Z}^d} p(e)p^d(e')P(\tau(0,e) = \tau(0,e') = \infty).$$

COROLLARY 1.10. $P_N \Rightarrow P_{X_0}^{2\gamma_e, \theta', \sigma^2}$ as $N \to \infty$, where $\theta' = \theta_0 \beta' - \theta_1 \delta'$.

PROOF. This is another application of Corollary 1.5 with

$$\beta_N(A) = \begin{cases} \theta_0^N p(a) p^b(a), & A = \{a\}, \\ \theta_0^N (p(a) p^b(a') + p(a') p^b(a)), & A = \{a, a'\}, \\ 0, & \text{otherwise,} \end{cases}$$

and

$$\delta_N(A) = \begin{cases} \theta_1^N, & A = \varnothing, \\ \theta_1^N (p(a) p^d(a) - p(a) - p^d(a)), & A = \{a\}, \\ \theta_1^N (p(a) p^d(a') + p(a') p^d(a)), & A = \{a, a'\}, \\ 0, & \text{otherwise.} \end{cases}$$

One proceeds by verifying the conditions of Corollary 1.6 and applying that result as in the proof of Theorem 1.2—the arguments are similar and left for the reader. $\square$

The outline of the rest of the paper is as follows. In Section 2 we derive some crude bounds on the size of $X_t^N(\mathbf{1})$, and obtain a semimartingale decomposition of $X_t^N(\phi)$ for a large class of test functions $\phi$. In Section 3 the proof of our main result is reduced to a moment bound (Proposition 3.3) and a key estimate (Proposition 3.4). Given these results, we establish tightness of our sequence $X_\cdot^N$, and show all limit points converge to super-Brownian motion with the given parameters. A comparison scheme with the biased voter model in Section 4 will give the above moment bound, and play an important role in the proof of the key estimate. The latter is proved in Section 6 after some necessary probability estimates are established in Section 5.

**2. Construction and decomposition.** Our goal in this section is to derive the martingale problem for $X_\cdot^N$ and derive some elementary bounds on $|\xi_t^N| = \sum_x \xi_t^N(x)$. We assume that $\xi_t^N$ is the spin-flip system with pregenerator $\Omega_N$ described in the previous section. In this section we will not need any of the kernel assumptions, and will only need (P5) and the following weaker form of (P1) of the perturbation assumptions:

(P1)''' $$\sum_{A \in P_F} (|\beta_N(A)| + |\delta_N(A)|) < \infty \qquad \text{for all } N.$$

Recall also that (H1) and (H2) hold as always. Throughout this section, $N$ will be fixed, and we will let $\mathcal{F}_t$ be the canonical right-continuous filtration associated with $\xi_t^N$. All martingales will be understood to be $\mathcal{F}_t$-martingales.

PROPOSITION 2.1.

$$E\left(\sup_{t\leq T}|\xi_t^N|^p\right) < \infty \qquad \text{for all } p > 0 \text{ and } T \in [0,\infty). \tag{2.1}$$

PROOF. Let $c_1 = \sum_{A \in P_F} |\beta_N(A)|$ [finite by (P1)$'''$], and let $\psi$ be a selection function on the nonempty subsets in $P_F$, that is, $\psi(A) \in A/\ell_N$ for all nonempty $A$. Define

$$\hat{c}(x,\eta) = N \sum_{e \in \mathbf{S_N}} p_N(e)\eta(x+e) + \sum_{A \in P_F} |\beta_N(A)|\eta(x+\psi(A)).$$

Let $\hat{\eta}(\cdot) \in \mathbb{Z}_+^{\mathbf{S_N}}$ be the pure birth particle system such that $\hat{\eta}(x) \to \hat{\eta}(x) + 1$ with rate $\hat{c}(x,\eta)$. Then $|\hat{\eta}_t| = \sum_x \hat{\eta}_t(x)$ is a pure birth process with birth rate $N + c_1$ for each particle (this makes the existence and uniqueness of this system starting from a configuration of finitely many ones obvious). If $\eta(x) = \mathbb{1}(\hat{\eta}(x) \geq 1)$, then $\eta$ is a spin-flip system with jump rate $c'(x,\eta) = \hat{c}(x,\eta)\mathbb{1}(\eta(x) = 0)$. It is easy to use (1.14) and (1.17) to see that if $\xi(x) = \eta(x) = 0$, then $c_N(x,\xi) \leq c'(x,\eta)$. If $\eta(x) = 1$, then $c_N(x,\xi) \geq 0 = c'(x,\eta)$. By Theorem III.1.5 of [11], if $\eta_0 = \xi_0^N$, we may construct versions of $\xi_\cdot^N$ and $\eta_\cdot$ so that with probability one, $\xi_t^N \leq \eta_t$ for all $t \geq 0$. [For $\xi, \xi' \in \{0,1\}^{\mathbf{S_N}}$, $\xi \leq \xi'$ means that $\xi(x) \leq \xi'(x)$ for all $x \in \mathbf{S_N}$.] This implies that

$$\sup_{t \leq T}|\xi_t^N| \leq \sup_{t \leq T}|\eta_t| = |\eta_T|.$$

(Here, it is easy to use (P1)$'''$ to check the condition (0.3) on page 122 of [11], and so Theorem III.1.5 may be applied.) Since the pure birth process $|\hat{\eta}_T|$ has moments of all orders (see, e.g., Example 6.8.4 in [8]), so does $|\eta_T|$ and the proof is complete. $\square$

PROPOSITION 2.2. *For all $x \in \mathbf{S_N}$ and $t \geq 0$,*

$$\xi_t^N(x) = \xi_0^N(x) + M_t^{N,x} + D_t^{N,x}, \tag{2.2}$$

*where $\{M_\cdot^{N,x}, x \in \mathbf{S_N}\}$ are orthogonal square-integrable martingales with predictable square functions given by*

$$\langle M^{N,x}\rangle_t = \int_0^t \bigg[\sum_{y \in \mathbf{S_N}} Np_N(y-x)(\xi_s^N(y) - \xi_s^N(x))^2$$

$$+ \sum_A \chi_N(A, x, \xi_s^N)(\beta_N(A)\mathbb{1}\{\xi_s^N(x) = 0\} \tag{2.3}$$

$$+ \delta_N(A)\mathbb{1}\{\xi_s^N(x) = 1\})\bigg]ds$$

*and*

$$D_t^{N,x} = \int_0^t \left[ \sum_{y \in \mathbf{S_N}} N p_N(y-x)(\xi_s^N(y) - \xi_s^N(x)) \right.$$

(2.4)
$$+ \sum_{A \in P_F} \chi_N(A, x, \xi_s^N)(\beta_N(A)\mathbb{1}\{\xi_s^N(x) = 0\}$$

$$\left. - \delta_N(A)\mathbb{1}\{\xi_s^N(x) = 1\}) \right] ds.$$

PROOF. We will use the fact (e.g., Theorem I.5.2 of [11]) that for $\phi$ in the domain of $\Omega_N$,

(2.5) $\quad M_t = \phi(\xi_t) - \phi(\xi_0) - \int_0^t \Omega_N \phi(\xi_s) \, ds \qquad$ is a martingale.

Letting $\phi_x(\xi) = \xi(x)$, a calculation shows that

$$\Omega_N \phi_x(\xi) = \sum_{y \in \mathbf{S_N}} N p_N(y - x)(\xi(y) - \xi(x))$$

$$+ \sum_{A \in P_F} \chi_N(A, x, \xi)[\beta_N(A)\mathbb{1}\{\xi(x) = 0\} - \delta_N(A)\mathbb{1}\{\xi(x) = 1\}].$$

An application of (2.5) now gives the decomposition in (2.2). It follows from (2.2) that $M_t^{N,x}$ is uniformly bounded on compact time intervals and, hence, square integrable.

To derive the facts about the square function, we proceed as follows. Define $\phi_{x,y}$ (in the domain of $\Omega_N$) by $\phi_{x,y}(\xi) = \xi(x)\xi(y)$, and apply Itô's formula to $\phi_{x,x}$. Since $(\xi_t^N(x))^2 = \xi_t^N(x)$, we obtain the (second) decomposition of $\xi_t^N(x)$,

$$\xi_t^N(x) = \xi_0^N(x) + 2\int_0^t \xi_{s-}^N(x) \, dD_s^{N,x} + 2\int_0^t \xi_{s-}^N(x) \, dM_s^{N,x} + [M^{N,x}]_t,$$

where $[M^{N,x}]_{\cdot}$ is the square variation function of $M_{\cdot}^{N,x}$. The stochastic integral above is a martingale, as is $[M^{N,x}]_t - \langle M^{N,x}\rangle_t$, and, hence,

$$\xi_t^N(x) - \xi_0^N(x) - 2\int_0^t \xi_{s-}^N(x) \, dD_s^{N,x} - \langle M^{N,x}\rangle_t$$

is a martingale. Thus, we have written $\xi_t^N(x)$ as the sum of a martingale and a continuous process of bounded variation in two ways. Equating the processes of bounded variation leads to

$$\langle M^{N,x}\rangle_t = D_t^{N,x} - 2\int_0^t \xi_{s-}^N(x) \, dD_s^{N,x}.$$

A short calculation now gives (2.3).

The proof that the martingales $M_t^{N,x}$ are orthogonal proceeds in the same way. We use (2.5) with $\phi = \phi_{x,y}$ to obtain a semimartingale decomposition for the product $\xi_t^N(x)\xi_t^N(y)$. We then apply Itô's formula to obtain a second decomposition. Equating the processes of bounded variation leads to $\langle M^{N,x}, M^{N,y}\rangle_t = 0$, and the proof is complete. □

With Proposition 2.2 in hand, we can now obtain a decomposition for $X_t^N(\phi)$. First we introduce the following notation. For

$$\psi \in C_b(\mathbf{S_N}), \qquad \phi = \phi_s(x), \qquad \dot{\phi}_s(x) \equiv \frac{\partial}{\partial s}\phi(s,x) \in C_b([0,T] \times \mathbf{S_N}),$$

and $s \leq T$, define

$$\mathcal{A}_N(\psi) = \sum_{y \in \mathbf{S_N}} Np_N(y-x)(\psi(y) - \psi(x)),$$

$$D_t^{N,1}(\phi) = \int_0^t X_s^N(\mathcal{A}_N\phi_s + \dot{\phi}_s)\,ds,$$

$$D_t^{N,2}(\phi) = \frac{1}{N}\int_0^t \sum_{x \in \mathbf{S_N}} \phi_s(x) \sum_{A \in P_F} \beta_N(A)\chi_N(A,x,\xi_s^N)\,ds,$$

$$D_t^{N,3}(\phi) = \frac{1}{N}\int_0^t \sum_{x \in \mathbf{S_N}} \phi_s(x) \sum_{A \in P_F} (\beta_N(A) + \delta_N(A))\xi_s^N(x)\chi_N(A,x,\xi_s^N)\,ds,$$

$$\langle M^N(\phi)\rangle_{1,t} = \frac{1}{N^2}\int_0^t \sum_{x \in \mathbf{S_N}} \phi_s^2(x) \sum_{y \in \mathbf{S_N}} Np_N(y-x)(\xi_s^N(y) - \xi_s^N(x))^2\,ds,$$

$$\langle M^N(\phi)\rangle_{2,t} = \frac{1}{N^2}\int_0^t \sum_{x \in \mathbf{S_N}} \phi_s^2(x) \sum_{A \in P_F} \chi_N(A,x,\xi_s^N)(\beta_N(A)\mathbb{1}\{\xi_s^N(x) = 0\}$$
$$+ \delta_N(A)\mathbb{1}\{\xi_s^N(x) = 1\})\,ds.$$

Note that $\langle M^N(\phi)\rangle_{2,t}$ may be negative.

PROPOSITION 2.3. *For $\phi, \dot{\phi} \in C_b([0,T] \times \mathbf{S_N})$ and $t \in [0,T]$,*

(2.6) $$X_t^N(\phi_t) = X_0^N(\phi_0) + D_t^N(\phi) + M_t^N(\phi),$$

*where*

(2.7) $$D_t^N(\phi) = D_t^{N,1}(\phi) + D_t^{N,2}(\phi) - D_t^{N,3}(\phi),$$

*and $M_t^N(\phi)$ is a square-integrable martingale with predictable square function*

(2.8) $$\langle M^N(\phi)\rangle_t = \langle M^N(\phi)\rangle_{1,t} + \langle M^N(\phi)\rangle_{2,t}.$$

PROOF. Use Proposition 2.2 and integration by parts to see that

$$\phi_t(x)\xi_t^N(x) = \phi_0(x)\xi_0^N(x) + \int_0^t \phi_s(x)\,dM_s^{N,x} + \int_0^t \phi_s(x)\,dD_s^{N,x}$$
(2.9)
$$+ \int_0^t \dot{\phi}_s(x)\xi_s^N(x)\,ds.$$

Using (P5) and the elementary inequality

(2.10) $$\chi_N(A, x, \xi_s^N) \leq \frac{1}{|A|} \sum_{a \in A/\ell_N} \xi_s^N(x+a), \qquad A \neq \varnothing,$$

we have

$$\sum_{x \in \mathbf{S_N}} \sum_{A \in P_F} \chi_N(A, x, \xi_s^N)(|\beta_N(A)|\mathbb{1}(\xi_s^N(x) = 0) + |\delta_N(A)|\mathbb{1}(\xi_s^N(x) = 1))$$

(2.11)
$$\leq \left[ \sum_{x \in \mathbf{S_N}} \sum_{A \in P_F} |A|^{-1} \sum_{a \in A} \xi_s^N(x + a/\ell_N)|\beta_N(A)| \right]$$
$$+ \left[ |\xi_s^N| \sum_{A \in P_F} |\delta_N(A)| \right]$$
$$\leq |\xi_s^N| \sum_{A \in P_F} (|\beta_N(A)| + |\delta_N(A)|).$$

This, together with Proposition 2.1, (H2) and (P1)′′′, shows that each of the terms in (2.9) is nonzero for only finitely many values of $x$ for all $t \leq T$ a.s. Here we first make this conclusion for each of the terms other than the martingale integral and, hence, infer it for the martingale integrals. We therefore may sum (2.9) over $x$, and after a bit of rearranging, obtain the required decomposition with

(2.12) $$M_t^N(\phi) = \frac{1}{N} \sum_{x \in \mathbf{S_N}} \int_0^t \phi_s(x)\,dM_s^{N,x}.$$

Now use (2.11) and Proposition 2.1 to see that

$$E\left( \sum_{x \in \mathbf{S_N}} \left\langle \int_0^\cdot \phi_s(x)\,dM_s^{N,x} \right\rangle_T \right) < \infty.$$

This shows that the series in (2.12) converges in $L^2$ uniformly in $t \leq T$ and so $M_N(\phi)$ is a square integrable martingale. It also shows that its predictable square function is

$$\lim_{K \to \infty} \frac{1}{N^2} \sum_{\substack{x \in \mathbf{S_N} \\ |x| \leq K}} \left\langle \int_0^\cdot \phi_s(x)\,dM_s^{N,x} \right\rangle_t,$$

where the limit exists in $L^1$ by the above but also for all $t \leq T$ a.s. by monotonicity. A simple calculation using (2.3) now gives (2.8) and the proof is complete. □

**3. Convergence to super-Brownian motion.** Our strategy in proving Theorem 1.3 is standard. We will prove that the family $\{X^N_\cdot, N \geq 1\}$ is tight, and that all weak limit points $X_\cdot$ satisfy the martingale problem characterizing super-Brownian motion $X_\cdot$ with the specified parameters. Hence, $X^N_\cdot \Rightarrow X_\cdot$ as $N \to \infty$. Our task here is less complicated than in [3], because we consider only the high-dimensional case, $d \geq 3$. The appropriate mass normalizer is $N' = N$, which fits well with Brownian space-time scaling. Many of the complications in [3] arose considering the delicate $d = 2$ case, for which the appropriate mass normalizer was $N' = N/\log N$. On the other hand, our task here is more difficult than in [3] because the Lotka–Volterra and perturbed voter models do not have tractable dual processes, as does the basic voter model.

A sequence of probability measures $\{P_N\}$ on $D([0,\infty), E)$ ($E$ a Polish space) is $C$-tight iff it is tight and every limit point is supported by $C([0,\infty), E)$. Recall that $P_N$ is the law of $X^N_\cdot$ on $D([0,\infty), M_f(\mathbb{R}^d))$, and that the assumptions of Theorem 1.3 are in force. Our strategy requires proving the following two results.

PROPOSITION 3.1. *The family of laws $\{P_N, N \in \mathbb{N}\}$ is $C$-tight.*

PROPOSITION 3.2. *If $P^*$ is any weak limit point of the sequence $P_N$, then $P^* = P^{2\gamma,\theta,\sigma^2}$.*

Clearly, Theorem 1.3 follows from these propositions.

We now state a pair of key technical results, Propositions 3.3 and 3.4 below, whose proofs we defer to Sections 4–6. Assuming these two propositions, we give the proofs of Propositions 3.1 and 3.2 in this section.

PROPOSITION 3.3. *For $K, T > 0$, there exists a finite constant $C_3(K,T)$ such that if $\sup_N X^N_0(\mathbf{1}) \leq K$, then*

$$\sup_N E\left(\sup_{t \leq T} X^N_t(\mathbf{1})^2\right) \leq C_3(K,T). \tag{3.1}$$

This bound allows us to employ $L^2$ arguments. Note that it is a consequence of (H2) that there will exist a $K$ as above.

Our second (and key) technical bound will need the following notation. For $A \in P_F$, $\phi : [0,T] \times \mathbf{S_N} \to \mathbb{R}$ bounded and measurable, $K > 0$ and $t \in [0,T]$,

define

$$\mathcal{E}_N(A, \phi, K, t)$$
$$= \sup_{X_0^N(\mathbf{1}) \leq K} E\left(\left(\int_0^t \left[\frac{1}{N}\sum_x \phi_s(x)\chi_N(A, x, \xi_s^N) - \sigma_N(A)X_s^N(\phi_s)\right] ds\right)^2\right)$$

[recall that $\sigma_N(A) = P(\tau_N(A) \leq \varepsilon_N^*)$]. For $\phi: \mathbf{S_N} \to \mathbb{R}$, define

$$\|\phi\|_{\text{Lip}} = \|\phi\|_\infty + \sup_{x \neq y} |\phi(x) - \phi(y)||x - y|^{-1}.$$

Also, recall that $\ell_N = M_N\sqrt{N} \to \infty$. By (P1), $c_\beta = \sup_N \sum_{A \in P_F} \beta_N(A)^+ < \infty$ and we may set $\bar{c} = c_\beta + k_\delta$, where $k_\delta$ is as in (P4).

PROPOSITION 3.4. *There is a positive sequence $\varepsilon_N \to 0$ as $N \to \infty$, and, for any $K, T > 0$, a constant $C_4(K, T) > 0$, such that for any $\phi \in \mathbf{C}_b([0, T] \times \mathbf{S_N})$ satisfying $\sup_{s \leq T} \|\phi_s\|_{\text{Lip}} \leq K$, nonempty $A \in P_F$, $\bar{a} \in A$, $J \geq 1$, and $0 \leq t \leq T$,*

$$\mathcal{E}_N(A, \phi, K, t) \leq C_4(K, T)[\varepsilon_N^* e^{\bar{c}\varepsilon_N^*} + J^{-2}$$
(3.2)
$$+ J^2(\varepsilon_N|A| + (\sigma_N(A) \wedge (\varepsilon_N + |\bar{a}|/\ell_N)))].$$

*In particular, $\lim_{N \to \infty} \sup_{t \leq T} \mathcal{E}_N(A, \phi, K, t) = 0$.*

This result says that

$$\frac{1}{N}\sum_x \phi_s(x)\chi_N(A, x, \xi_s^N) \approx \sigma_N(A)X_s^N(\phi_s),$$

in some average sense, and is the key to identifying any weak limit of $X_\cdot^N$. We proceed now assuming the validity of the above two propositions.

We begin by obtaining more precise information on the terms in the decomposition of $X_t^N(\phi)$ given in Proposition 2.3. Lemma 3.5 below estimates the terms in the increasing process $\langle M^N(\phi)\rangle_t$, Lemma 3.6 estimates the terms in the drift $D_t^N(\phi)$.

LEMMA 3.5. *There is a constant $C$ such that if $\phi: [0, T] \times \mathbf{S_N} \to \mathbb{R}$ is a bounded measurable function, then*

(a) $\langle M^N(\phi)\rangle_{2,t} = \int_0^t m_{2,s}^N(\phi)\,ds$, *where*

(3.3) $$|m_{2,s}^N(\phi)| \leq C\frac{\|\phi\|_\infty^2}{N}X_s^N(\mathbf{1}).$$

(b)

$$\langle M^N(\phi)\rangle_{1,t} = 2\int_0^t X_s^N(\phi_s^2 f_0^N(\xi_s^N))\,ds + \int_0^t m_{1,s}^N(\phi_s)\,ds, \tag{3.4}$$

where

$$|m_{1,s}^N(\phi)| \leq \left[\frac{C}{\sqrt{N}}\|\phi_s\|_{\text{Lip}}^2 X_s^N(\mathbf{1})\right] \wedge [2\|\phi\|_\infty^2 X_s^N(\mathbf{1})]. \tag{3.5}$$

(c) For $i=2,3$, $D_t^{N,i}(\phi) = \int_0^t d_s^{N,i}(\phi)\,ds$ for $t \leq T$, where for all $N$ and $s \leq T$,

$$|d_s^{N,i}(\phi)| \leq C\|\phi\|_\infty X_s^N(\mathbf{1}).$$

PROOF. (a) The definition of $\langle M^N(\phi)\rangle_{2,t}$ implies

$$|m_{2,s}^N(\phi)| \leq \frac{1}{N^2}\sum_{x\in\mathbf{S_N}}|\phi_s(x)|^2 \sum_{A\in P_F\setminus\varnothing}(|\beta_N(A)|+|\delta_N(A)|)\chi_N(A,x,\xi_s^N)$$

$$+ \frac{1}{N}X_s^N(\phi_s^2)|\delta_N(\varnothing)|.$$

By (P1) and (2.10), there is a constant $C$ such that

$$|m_{2,s}^N(\phi)| \leq \|\phi\|_\infty^2 \sum_{A\in P_F\setminus\varnothing}\frac{(|\beta_N(A)|+|\delta_N(A)|)}{|A|}\sum_{a\in A}\frac{1}{N^2}\sum_{x\in\mathbf{S_N}}\xi_s^N\left(x+\frac{a}{\ell_N}\right)$$

$$+ \|\phi\|_\infty^2 \frac{\delta_N(\varnothing)}{N}X_s^N(\mathbf{1}) \tag{3.6}$$

$$\leq C\frac{\|\phi\|_\infty^2}{N}X_s^N(\mathbf{1}).$$

(c) This is proved by making minor changes in the derivation of (a).

(b) A little rearrangement is necessary to handle the term $\langle M^N(\phi)\rangle_{1,t}$. We rewrite it in the form

$$\frac{1}{N}\int_0^t \sum_{x,y\in\mathbf{S_N}} p_N(y-x)\phi_s^2(x)[\xi_s^N(x)(1-\xi_s^N(y))+(1-\xi_s^N(x))\xi_s^N(y)]\,ds$$

$$= \frac{1}{N}\int_0^t \sum_{x,y\in\mathbf{S_N}} \xi_s^N(x)\phi_s^2(x)p_N(y-x)(1-\xi_s^N(y))\,ds$$

$$+ \frac{1}{N}\int_0^t \sum_{x,y\in\mathbf{S_N}} \xi_s^N(y)\phi_s^2(y)p_N(y-x)(1-\xi_s^N(x))\,ds$$

$$+ \frac{1}{N}\int_0^t \sum_{x,y\in\mathbf{S_N}} p_N(y-x)[\phi_s^2(x)-\phi_s^2(y)]\xi_s^N(y)(1-\xi_s^N(x))\,ds.$$

That is, (3.4) holds where

$$m_{1,s}^N(\phi) = \frac{1}{N} \sum_{x,y \in \mathbf{S_N}} p_N(y-x)(\phi_s^2(x) - \phi_s^2(y))\xi_s^N(y)(1 - \xi_s^N(x)).$$

Note that $|\phi_s(x)^2 - \phi_s(y)^2| \leq 2\|\phi_s\|_{\text{Lip}}^2 |x-y|$, and also, by (H1) for some universal constant $C$,

$$\sum_y p_N(y-x)|x-y| = E(|W_N|)/\sqrt{N} \leq C/(2\sqrt{N}).$$

These inequalities establish (3.5). □

Let $T > 0$ and $\phi : [0,T] \times \mathbf{S_N} \to \mathbb{R}$ be such that $\phi, \dot{\phi} \in C_b([0,T] \times \mathbf{S_N})$, and define

$$\delta_N^1(s,\phi) = \sum_{A \in P_F} \beta_N(A) \left[ \frac{1}{N} \sum_{x \in \mathbf{S_N}} \phi_s(x)\chi_N(A, x, \xi_s^N) - \sigma_N(A)X_s^N(\phi_s) \right],$$

$$\delta_N^2(s,\phi) = \sum_A (\beta_N(A) + \delta_N(A)) \left[ \frac{1}{N} \sum_{x \in \mathbf{S_N}} \phi_s(x)\chi_N(A \cup \{0\}, x, \xi_s^N) \right.$$
$$\left. - \sigma_N(A \cup \{0\})X_s^N(\phi_s) \right].$$

It follows from (2.10), (P1), (P5) and Proposition 3.3 that these series converge. Also, set

$$d_0^N = \sum_{A \in P_F} \beta_N(A)\sigma_N(A) - \sum_{A \in P_F} (\beta_N(A) + \delta_N(A))\sigma_N(A \cup \{0\}),$$

and note by (P1) that

(3.7) $$c_1 = \sup_N |d_0^N| < \infty.$$

With this notation, (2.6) of Proposition 2.3 may be written as

$$X_t^N(\phi_t) = X_0^N(\phi_0) + M_t^N(\phi) + \int_0^t X_s^N(\mathcal{A}_N \phi_s + \dot{\phi}_s)\, ds$$

(3.8) $$+ \int_0^t d_0^N X_s^N(\phi_s)\, ds + \int_0^t (\delta_N^1(s,\phi) - \delta_N^2(s,\phi))\, ds$$

$$\text{for all } t \in [0,T].$$

LEMMA 3.6. *There is a sequence $\varepsilon_N^0 \to 0$ as $N \to \infty$ and for each $K, T > 0$ a constant $C_0(K,T)$ (increasing in each variable) such that if $\phi: [0,T] \times \mathbf{S_N} \to \mathbb{R}$ satisfies $\sup_{s \leq T} \|\phi_s\|_{\text{Lip}} \leq K$ and $\sup_N X_0^N(\mathbf{1}) \leq K$, then*

$$(3.9) \quad \sup_{t \leq T}\left[E\left(\left(\int_0^t \delta_N^1(s,\phi)\,ds\right)^2 + \left(\int_0^t \delta_N^2(s,\phi)\,ds\right)^2\right)\right]^{1/2} \leq C_0(T,K)\varepsilon_N^0$$

*for all $N$.*

PROOF. Assume $\phi$ and $X_0^N$ are as above. If $t \in [0,T]$, then by Cauchy–Schwarz and (P1),

$$\begin{aligned}
E\left(\left(\int_0^t \delta_N^1(s,\phi)\,ds\right)^2\right) \\
= E\left(\left(\sum_{A \in P_F} \beta_N(A) \int_0^t \left[\frac{1}{N}\sum_{x \in \mathbf{S_N}} \phi_s(x)\chi_N(A,x,\xi_s^N) \right.\right.\right. \\
\left.\left.\left. - \sigma_N(A)X_s^N(\phi_s)\right]ds\right)^2\right) \\
\leq C \sum_{A \in P_F} |\beta_N(A)| \mathcal{E}_N(A,\phi,K,t)
\end{aligned} \tag{3.10}$$

for a constant $C$. Proposition 3.4 and (P1) show that for some positive sequence $\varepsilon_N' \to 0$ and any $J \geq 1$,

$$\sup_{t \leq T} E\left(\left(\int_0^t \delta_N^1(s,\phi)\,ds\right)^2\right) \leq C(T,K)(\varepsilon_N' + J^{-2} + J^2(\varepsilon_N' + \eta_N)),$$

where $C(T,K)$ does not depend on the choice of $\phi$, and

$$\eta_N = \sum_A |\beta_N(A)|(\sigma_N(A) \wedge (\varepsilon_N + |\bar{a}|/\ell_N)).$$

(Recall $\bar{a}$ denotes some element of $A$.) By (P3) and a uniform integrability argument, $\eta_N \to 0$ as $N \to \infty$. Optimize the above over $J$ to see that for some positive sequence $\varepsilon_N'' \to 0$,

$$\sup_{t \leq T} E\left(\left(\int_0^t \delta_N^1(s,\phi)\,ds\right)^2\right) \leq C(T,K)\varepsilon_N''.$$

A similar argument goes through for $\delta_N^2(s,\phi)$ [note that $\sigma_N(A \cup \{0\}) \leq \sigma_N(A)$] and so the result follows (the monotonicity requirements on $C_0$ are trivial to realize). □

The proof of Proposition 3.1 (tightness) proceeds as follows. We first establish tightness for $X_\cdot^N(\phi)$ for an appropriate class of test functions $\phi$. We then prove a "compact containment" condition for $X_\cdot^N$. We can then appeal to a version of Jakubowski's theorem for weak convergence in $D([0,\infty), \mathcal{M}_f(\mathbb{R}^d))$ (see Theorem II.4.1 in [14]), completing the proof of Proposition 3.1.

PROPOSITION 3.7. *For each* $\phi \in C_b^{1,3}(\mathbb{R}_+ \times \mathbb{R}^3)$, *each of the families* $\{X_\cdot^N(\phi_\cdot), N \in \mathbb{N}\}$, $\{D_\cdot^N(\phi), N \in \mathbb{N}\}$, $\{\langle M^N(\phi)\rangle_\cdot, N \in \mathbb{N}\}$ *and* $\{M_\cdot^N(\phi), N \in \mathbb{N}\}$ *is $C$-tight in* $D([0,\infty), \mathbb{R})$.

PROOF. Fix $\phi$ as above and recall the decomposition of $X_t^N(\phi_t)$ in Proposition 2.3. We start with the drift terms and recall an analytic estimate (Lemma 2.6) of [3]:

$$(3.11) \qquad \sup_{s \leq T} \left\| \mathcal{A}_N(\phi_s) - \frac{\sigma^2 \Delta \phi_s}{2} \right\|_\infty \to 0 \qquad \text{as } N \to \infty.$$

Since $D_t^{N,1}(\phi) = \int_0^t X_s^N(\mathcal{A}_N \phi_s + \dot\phi)\, ds$, (3.11), Proposition 3.3 and the Arzela–Ascoli theorem imply that

$$\{D_\cdot^{N,1}(\phi), N \in \mathbb{N}\} \qquad \text{is tight in } C([0,\infty), \mathbb{R}).$$

For $i = 2, 3$, $D_t^{N,i}(\phi) = \int_0^t d_s^{N,i}(\phi)\, ds$, where by Lemma 3.5(c),

$$|d_s^{N,i}(\phi)| \leq C\|\phi\|_\infty X_s^N(\mathbf{1}), \qquad i = 2, 3.$$

Again Proposition 3.3 and the Arzela–Ascoli theorem imply that

$$\{D_\cdot^{N,i}(\phi), N \in \mathbb{N}\} \qquad \text{is tight in } C([0,\infty), \mathbb{R}), i = 2, 3.$$

We turn now to the martingale terms. By (2.8) and Lemma 3.5(a, b), there is a finite constant $C$ such that for $0 \leq s \leq t \leq T$,

$$(3.12) \qquad \langle M^N(\phi)\rangle_t - \langle M^N(\phi)\rangle_s \leq C\|\phi\|_\infty^2 \int_s^t X_u^N(\mathbf{1})\, du.$$

Consequently, Proposition 3.3 shows that

$$\{\langle M^N(\phi)\rangle_\cdot, N \in \mathbb{N}\} \qquad \text{is tight in } C([0,\infty), \mathbb{R}).$$

Since the maximum jump discontinuity in $M_t^N(\phi)$ is bounded above by $\|\phi\|_\infty/N$, it follows from Theorem VI.4.13 and Proposition VI.3.26 of [10] that

$$\{M_\cdot^N(\phi), N \in \mathbb{N}\} \qquad \text{is $C$-tight in } D([0,\infty), \mathbb{R}).$$

In view of (H2), we see from the above and Proposition 2.3 that $X_t^N(\phi_t)$ and $D_t^N(\phi)$ are each a sum of $C$-tight processes in $D([0,\infty), \mathbb{R})$. Since a sum of $C$-tight processes in $D([0,\infty), \mathbb{R})$ is also $C$-tight, the proof is complete. □

To derive the appropriate compact containment condition, we will first need an estimate on the mean measure of $X_t^N$. Let $P_t^N$ denote the semigroup associated with the generator $\mathcal{A}_N$.

PROPOSITION 3.8. *There is a constant $c_1 \geq 0$, a positive sequence $\epsilon_N^1 \to 0$ as $N \to \infty$, and constants $(C_1(K,t), K, t \geq 0)$, nondecreasing in each variable, such that if $\sup_N X_0^N(\mathbf{1}) \leq K$, and $\phi : \mathbf{S_N} \to \mathbb{R}_+$ satisfies $\|\phi\|_{\mathrm{Lip}} \leq K$, then*
$$E(X_t^N(\phi)) \leq e^{c_1 t} X_0^N(P_t^N \phi) + C_1(K,t)\epsilon_N^1.$$

PROOF. Assume $c_1$ is as in (3.7) and $\phi$ is as in the statement of the proposition. Fix $t > 0$ and define
$$\phi_s(x) = e^{-c_1 s} P_{t-s}^N \phi(x), \qquad (s,x) \in [0,t] \times \mathbf{S_N}.$$
Then (3.8) becomes
$$e^{-c_1 t} X_t^N(\phi) = X_0^N(P_t^N \phi) + M_t^N(\phi) + (d_0^N - c_1) \int_0^t X_s^N(\phi_s)\,ds$$
$$+ \int_0^t (\delta_N^1(s,\phi) - \delta_N^2(s,\phi))\,ds.$$
Note that the third term on the right-hand side is nonpositive. It is easy to verify that $\sup_{s \leq t} \|\phi_s\|_{\mathrm{Lip}} \leq K$. Therefore, we may use Lemma 3.6, and take expectations in the above with $T = t$, recalling that $M_t^N(\phi)$ is a mean zero martingale (Proposition 2.3), to arrive at
$$E(X_t^N(\phi)) \leq e^{c_1 t} X_0^N(P_t^N \phi) + e^{c_1 t} 2 C_0(K,t)\epsilon_N^0.$$
The result is then immediate. $\square$

For the following, let $B(x,r)$ denote the open ball in $\mathbb{R}^d$ of radius $r$ centered at $x$.

PROPOSITION 3.9 (Compact containment). *For all $\epsilon > 0$, there is a finite $\rho = \rho(\epsilon)$ such that*
$$\sup_N P\left(\sup_{t \leq \epsilon^{-1}} X_t^N(B(0,\rho)^c) > \epsilon\right) < \epsilon.$$

PROOF. Let $h_n : \mathbb{R}^d \to [0,1]$ be a $C^\infty$ function such that
$$B(0,n) \subset \{x : h_n(x) = 0\} \subset \{x : h_n(x) < 1\} \subset \{B(0,n+1)\}$$
and
$$\sup_n \sum_{i,j,k \leq d} \|(h_n)_i\|_\infty + \|(h_n)_{ij}\|_\infty + \|(h_n)_{ijk}\|_\infty \equiv C_h < \infty.$$

Let $c_1$ be as in (3.7) and use (3.8) with $\phi_s^n(x) = e^{-c_1 s} h_n(x)$ to get

$$
\begin{aligned}
(3.13) \quad e^{-c_1 t} X_t^N(h_n) = {} & X_0^N(h_n) + M_t^N(\phi^n) + \int_0^t e^{-c_1 s} X_s^N(\mathcal{A}_N h_n) \\
& + (d_0^N - c_1) X_s^N(\phi_s^n) \, ds + \int_0^t \delta_N^1(s, \phi^n) - \delta_N^2(s, \phi^n) \, ds.
\end{aligned}
$$

Note that

$$
\begin{aligned}
(3.14) \quad E\!\left(\int_0^t X_s^N(|\mathcal{A}_N h_n|) \, ds\right) \leq {} & \left\| \mathcal{A}_N h_n - \frac{\sigma^2 \Delta h_n}{2} \right\|_\infty E\!\left(\int_0^t X_s^N(\mathbf{1}) \, ds\right) \\
& + E\!\left(\int_0^t X_s^N\!\left(\frac{\sigma^2 |\Delta h_n|}{2}\right) ds\right).
\end{aligned}
$$

The first term in (3.14) approaches zero as $N \to \infty$, uniformly in $n$ by (3.11) and Proposition 3.3. Choose

$$
(3.15) \quad K > \max\!\left(1, C_h(\sigma^2/2 + 1), \sup_N X_0^N(\mathbf{1})\right).
$$

Then $\phi = \sigma^2 |\Delta h_n|/2$ satisfies the hypotheses of Proposition 3.8 and so that result bounds the second term in (3.14) by

$$
(3.16) \quad \int_0^t e^{c_1 s} X_0^N\!\left(P_s^N\!\left(\frac{\sigma^2 |\Delta h_n|}{2}\right)\right) ds + C_1(K, t) t \epsilon_N^1.
$$

Since $\Delta h_n = 0$ on $B(0, n)$, we may use (H1) and (H2) to conclude that

$$
\begin{aligned}
X_0^N(P_s^N(|\Delta h_n|)) &\leq C_h X_0^N(P_s^N(\mathbf{1}_{B(0,n)^c})) \\
&\leq C_h(X_0^N(B(0, n/2)^c) + X_0^N(\mathbf{1}) P(|B_s^{0,N}| > n/2)) \\
&\leq C_h(X_0^N(B(0, n/2)^c) + X_0^N(\mathbf{1}) c n^{-2} s) \\
&\to 0 \quad \text{as } n \to \infty \text{ uniformly in } N \text{ and } s \leq t.
\end{aligned}
$$

The above proves

$$
(3.17) \quad \lim_{(N,n) \to \infty} E\!\left(\int_0^t X_s^N(|\mathcal{A}_N h_n|) \, ds\right) = 0.
$$

Use (2.8) and Lemma 3.5 to see that [recall $\phi_s^n(x) = e^{-c_1 s} h_n(x)$]

$$
\begin{aligned}
(3.18) \quad E(\langle M^N(\phi^n)\rangle_t) \leq {} & C(N^{-1} + N^{-1/2}) E\!\left(\int_0^t X_s^N(\mathbf{1}) \, ds\right) \\
& + 2 E\!\left(\int_0^t X_s^N(h_n^2) \, ds\right).
\end{aligned}
$$

Now use Proposition 3.8 to bound the second term in (3.18) [just as in (3.16)] and Proposition 3.3 to bound the first term in (3.18) and conclude

$$\lim_{(N,n)\to\infty} E(\langle M^N(\phi^n)\rangle_t) = 0 \quad \text{for all } t > 0. \tag{3.19}$$

Let $\epsilon > 0$. By (H2), (3.17) and (3.19) there is an $n_0 \in \mathbb{N}$ such that for $N, n \geq n_0$,

$$P\left(e^{c_1\epsilon^{-1}} X_0^N(h_n) + \sup_{t \leq \epsilon^{-1}} e^{c_1 t} |M_t^N(\phi^n)| \right.$$
$$\left. + \int_0^{\epsilon^{-1}} e^{c_1(t-s)} X_s^N(|\mathcal{A}_N h_n|) \, ds > \epsilon \right) < \epsilon. \tag{3.20}$$

Turning now to the last term in (3.13), note first the trivial bound

$$|\delta_N^1(s, \phi^n)| + |\delta_N^2(s, \phi^n)| \leq \sum_A (|\beta_N(A)| + |\delta_N(A)|) 4 X_s^N(\mathbf{1})$$
$$\leq C X_s^N(\mathbf{1}), \tag{3.21}$$

the last inequality by (P1). Our choice of $K$ in (3.15) shows that each $\phi^n$ satisfies $\sup_s \|\phi_s^n\|_{\text{Lip}} \leq K$ and so Lemma 3.6 implies that for all $T > 0$,

$$\sup_{t \leq T} E\left(\left|\int_0^t \delta_N^i(s, \phi^n) \, ds\right|\right) \to 0 \tag{3.22}$$

as $N \to \infty$ uniformly in $n$ for $i = 1, 2$.

Now (3.21) and Proposition 3.3 show that $\{\int_0^\cdot \delta_N^i(s, \phi^{n_0}) \, ds : N \in \mathbb{N}\}$, $i = 1, 2$, are tight in $C(\mathbb{R}_+, \mathbb{R})$, while (3.22) shows that each limit point of the above sequences is identically 0. This shows weak convergence of $\int_0^\cdot \delta_N^i(s, \phi^{n_0}) \, ds$ to the zero process and, therefore,

$$\lim_{N\to\infty} P\left(\sup_{t\leq \epsilon^{-1}} e^{c_1 t} \left\{\left|\int_0^t \delta_N^1(s, \phi^{n_0}) \, ds\right| + \left|\int_0^t \delta_N^2(s, \phi^{n_0}) \, ds\right|\right\} > \epsilon\right) = 0.$$

Now use the above and (3.20) in (3.13), noting that $(d_0^N - c_1) X_s^N(\phi_s^{n_0}) \leq 0$, and conclude that there is an $N_0$ so that if $N \geq N_0$,

$$P\left(\sup_{t \leq \epsilon^{-1}} X_t^N(h_{n_0}) > 2\epsilon\right) < 2\epsilon.$$

By increasing $n_0$ if necessary to handle $N \leq N_0$, we get

$$\sup_N P\left(\sup_{t \leq \epsilon^{-1}} X_t^N(h_{n_0}) > 2\epsilon\right) < 2\epsilon,$$

and the proof is complete because $h_{n_0} \geq \mathbb{1}_{B(0,n_0+1)^c}$. □

PROOF OF PROPOSITION 3.1. The $C$-tightness of $\{P_N, N \in \mathbb{N}\}$ is now immediate from Propositions 3.7 and 3.9 above, and Theorem II.4.1 in [14]. □

PROOF OF PROPOSITION 3.2. We assume below that $\phi \in C_b^{1,3}([0,T] \times \mathbb{R}^d)$, $\sup_N X_0^N(\mathbf{1}) \leq K$ [such a $K$ exists by (H2)] and $0 \leq t \leq T$. First, (3.11) and Proposition 3.3 imply

$$(3.23) \quad E\left(\left(D_t^{N,1}(\phi) - \int_0^t X_s^N\left(\frac{\sigma^2 \Delta \phi_s}{2} + \dot{\phi}_s\right) ds\right)^2\right) \to 0 \qquad \text{as } N \to \infty.$$

We also have

$$D_t^{N,2}(\phi) - D_t^{N,3}(\phi) = \int_0^t \delta_N^1(s,\phi) - \delta_N^2(s,\phi)\, ds + d_0^N \int_0^t X_s^N(\phi)\, ds.$$

It follows from (P3), $\sigma_N(A \cup \{0\}) \leq \sigma_N(A)$, (P2) and (K3) that

$$\beta_N(\cdot)\sigma_N(\cdot \cup \{0\}) \to \beta(\cdot)\sigma(\cdot \cup \{0\}) \qquad \text{in } \ell_1(P_F) \text{ as } N \to \infty.$$

This and (P3) imply that $d_0^N \to \theta$ as $N \to \infty$. We may apply these results with Proposition 3.3 and Lemma 3.6 to conclude

$$(3.24) \quad E\left(\left(D_t^{N,2}(\phi) - D_t^{N,3}(\phi) - \theta \int_0^t X_s^N(\phi_s)\, ds\right)^2\right) \to 0 \qquad \text{as } N \to \infty.$$

We claim now that

$$(3.25) \quad E\left(\left(\langle M^N(\phi)\rangle_t - 2\gamma \int_0^t X_s^N(\phi_s^2)\, ds\right)^2\right) \to 0 \qquad \text{as } N \to \infty.$$

Define

$$\gamma_N = \sum_{e \in \mathbf{S_N}} p_N(e) P(\hat{\tau}^N(\{0,e\}) > \varepsilon_N^*)$$

[recall $\tau^N(A) = \hat{\tau}^N(A/\ell_N)$ for $A \subset \mathbb{Z}^d$]. By (2.8), Lemma 3.5, (K3) and Proposition 3.3, to prove (3.25), it suffices to prove that

$$(3.26) \quad E\left(\left(\int_0^t X_s^N(\phi_s^2 f_1^N(\xi_s^N)) - (1 - \gamma_N) X_s^N(\phi_s^2)\, ds\right)^2\right) \to 0 \qquad \text{as } N \to \infty.$$

To do this, we expand the integrand above in the form
$X_s^N(\phi_s^2 f_1^N(\xi_s^N)) - (1 - \gamma_N) X_s^N(\phi_s^2)$

$$= \frac{1}{N} \sum_{x \in \mathbf{S_N}} \phi_s^2(x) \xi_s^N(x) \sum_{y \in \mathbf{S_N}} p_N(y-x)[\xi_s^N(y) - P(\tau^N(0,(y-x)\ell_N) \leq \varepsilon_N^*)]$$

$$= \sum_{y \in \mathbf{S_N}} p_N(y) \frac{1}{N} \sum_{x \in \mathbf{S_N}} \phi_s^2(x) \xi_s^N(x)[\xi_s^N(x+y) - P(\tau^N(0,y\ell_N) \leq \varepsilon_N^*)]$$

$$= \sum_{a \in \mathbb{Z}^d} p_N(a/\ell_N)\left[\frac{1}{N} \sum_{x \in \mathbf{S_N}} \phi_s^2(x) \chi_N(\{0,a\}, x, \xi_s^N) - \sigma_N(\{0,a\}) X_s^N(\phi_s^2)\right].$$

Applying Cauchy–Schwarz, the left-hand side of (3.26) is bounded above by

$$\sum_{a \in \mathbb{Z}^d} p_N(a/\ell_N) E\left(\left(\int_0^t \left[\frac{1}{N} \sum_{x \in \mathbf{S_N}} \phi_s^2(x) \chi_N(\{0, a\}, x, \xi_s^N) - \sigma_N(\{0, a\}) X_s^N(\phi_s^2)\right] ds\right)^2\right).$$

Proposition 3.4 now completes the proof of (3.26) and, hence, of (3.25).

The above $L^2$ estimates [i.e., (3.23)–(3.25)] imply that for $\varepsilon > 0$,

$$P\left(\left|D_t^N(\phi) - \int_0^t X_s^N\left(\frac{\sigma^2}{2}\Delta\phi + \dot\phi_s\right) ds - \theta \int_0^t X_s^N(\phi) \, ds\right| > \varepsilon\right) \to 0$$

and

$$P\left(\left|\langle M^N(\phi)\rangle_t - 2\gamma \int_0^t X_s^N(\phi^2) \, ds\right| > \varepsilon\right) \to 0$$

as $N \to \infty$.

Now suppose that $P(X_\cdot^{N_k} \in \cdot) \Rightarrow P(X_\cdot \in \cdot)$ in $D([0,\infty), \mathcal{M}_f(\mathbb{R}^d))$ for some $X_\cdot \in C([0,\infty), \mathcal{M}_f(\mathbb{R}^d))$ as $k \to \infty$. Since $(X_\cdot^{N_k}, D_\cdot^{N_k}(\phi), \langle M^{N_k}(\phi)\rangle_\cdot)$ is $C$-tight in $D([0,\infty), \mathcal{M}_F(\mathbb{R}^d) \times C(\mathbb{R}) \times C(\mathbb{R}^+))$ [by Theorem 3.7 and Proposition 3.1], by Skorohod's theorem (taking a further subsequence if necessary), we may assume that

$$(X_\cdot^{N_k}, D_\cdot^{N_k}(\phi), \langle M_\cdot^{N_k}(\phi)\rangle) \to (X_\cdot, D_\cdot(\phi), L_\cdot(\phi)) \qquad \text{a.s.},$$

where $(X_\cdot, D_\cdot(\phi), L_\cdot(\phi))$ is continuous. By the probability estimates above, it follows that

$$(3.27) \quad D_t(\phi) = \int_0^t X_s\left(\frac{\sigma^2}{2}\Delta\phi + \dot\phi_s\right) ds + \theta \int_0^t X_s(\phi_s) \, ds \qquad \forall t \geq 0 \text{ a.s.}$$

and

$$(3.28) \qquad L_t(\phi) = 2\gamma \int_0^t X_s(\phi_s^2) \, ds \qquad \forall t \geq 0 \text{ a.s.}$$

By Proposition 2.3, $M_\cdot^{N_k}(\phi) \to M_\cdot(\phi) \in C(\mathbb{R})$ a.s., where

$$X_t(\phi_t) = X_0(\phi_0) + M_t(\phi) + \int_0^t X_s\left(\frac{\sigma^2 \Delta\phi_s}{2} + \dot\phi_s\right) ds + \int_0^t X_s(\theta\phi_s) \, ds,$$

(3.29)

and $M_t(\phi)$ is continuous and $\mathcal{F}_t^X$-measurable. By (3.25) and Proposition 3.3,

$$\sup_N E(\langle M^N(\phi)\rangle_T^2) < \infty.$$

Using Burkholder's inequality and the fact that $|\Delta M^N(\phi)(t)| \leq \|\phi\|_\infty/N$, we obtain

$$\sup_N E\left(\sup_{t \leq T} |M_t^N(\phi)|^4\right) < \infty.$$

Consequently, $M.(\phi)$ is a continuous, $L^2$, $\mathcal{F}_\cdot^X$-measurable martingale, and

$$\langle M(\phi)\rangle_t = \lim_{k \to \infty} \langle M^{N_k}(\phi)\rangle_t = 2\gamma \int_0^t X_s(\phi_s^2)\, ds \qquad \text{a.s.}$$

Consequently, $P(X. \in \cdot)$ satisfies the martingale problem characterizing $P^{2\gamma,\theta,\sigma^2}$, and so $P(X_\cdot^{N_k} \in \cdot) \Rightarrow P^{\gamma,\theta,\sigma^2}$ as $N_k \to \infty$. $\square$

**4. Comparison with biased voter models.** In this section we show that we can dominate the process $\xi_t^N$ by a biased voter model $\bar{\xi}_t^N$. That is, we show that the two processes can be coupled so that with probability one, $\xi_t^N \leq \bar{\xi}_t^N$ for all $t \geq 0$. Easily obtained bounds on $E(|\bar{\xi}_t^N|)$ and $E(|\bar{\xi}_t^N|^2)$ thus provide bounds on $E(X_t^N(\mathbf{1}))$ and $(E(X_t^N(\mathbf{1})))^2$. The results in this section will use (P1), (P4) and (P5), but not any of the kernel assumptions.

Let $p$ and $\bar{p}$ be two probability kernels on $\mathbb{Z}^d$, and fix parameters $v > 0, b \geq 0$. For $i = 0, 1$, define

$$f_i(x, \eta) = \sum_{y \in \mathbb{Z}^d} p(y - x)\mathbb{1}\{\eta(y) = i\}$$

and

$$\bar{f}_i(x, \eta) = \sum_{y \in \mathbb{Z}^d} \bar{p}(y - x)\mathbb{1}\{\eta(y) = i\}.$$

The biased voter model $\bar{\xi}_t$ is the spin-flip system taking values in $\{0,1\}^{\mathbb{Z}^d}$ which in state $\bar{\xi}$ makes transitions at $x$,

(4.1)
$$\begin{aligned} 0 \to 1 &\quad \text{at rate } vf_1(x, \bar{\xi}) + b\bar{f}_1(x, \bar{\xi}), \\ 1 \to 0 &\quad \text{at rate } vf_0(x, \bar{\xi}). \end{aligned}$$

If $b = 0$, we obtain the voter model, while if $b > 0$, there is a bias in favor of creating 1's. It is clear from these rates that we may as well assume $p(0) = \bar{p}(0) = 0$.

We will need the following estimates on the first two moments of $|\bar{\xi}_t|$.

LEMMA 4.1. *Assume that $|\bar{\xi}_0| < \infty$. Then $|\bar{\xi}_t|$ is submartingale such that*

(4.2) $$E(|\bar{\xi}_t|) \leq e^{bt}|\bar{\xi}_0|,$$

*and $|\bar{\xi}_t|^2$ is a submartingale such that*

(4.3) $$E(|\bar{\xi}_t|^2) \leq e^{2bt}\left(|\bar{\xi}_0|^2 + \frac{b + 2v}{b}(1 - e^{-bt})|\bar{\xi}_0|\right).$$

PROOF. First, note that by bounding $|\bar{\xi}_t|$ above by a pure birth process just as in the proof of Proposition 2.1, one may conclude that for $T > 0$, the first and second moments of $\sup_{t \leq T} |\bar{\xi}_t|$ are finite. Next, if $\beta_1(\{a\}) = \frac{b}{v}\bar{p}(a)$, $\beta_1(A) = 0$ if $|A| \neq 1$, and $\delta_1 \equiv 0$, then $|\bar{\xi}_{t/v}|$ is precisely $X_t^1(\mathbf{1})$, where $X_\cdot^1$ is as in Theorem 1.3 with $N = 1$. Clearly, $\beta_1(A) = 0$ if $0 \in A$, (P5) holds and (P1)''' is valid, so from Proposition 2.3,

$$(4.4) \qquad |\bar{\xi}_t| = |\bar{\xi}_0| + \int_0^t \sum_{x,e \in \mathbb{Z}^d} b\bar{p}(e) \bar{\xi}_s(x+e)(1 - \bar{\xi}_s(x)) \, ds + \bar{M}_t,$$

where $\bar{M}_t$ is a square-integrable martingale with predictable square function

$$(4.5) \qquad \langle \bar{M} \rangle_t = \int_0^t \sum_{x,y \in \mathbb{Z}^d} \left[ vp(y-x) \mathbb{1}(\bar{\xi}_s(x) \neq \bar{\xi}_s(y)) \right.$$
$$\left. + \sum_{x,e \in \mathbb{Z}^d} b\bar{p}(e) \bar{\xi}_s(x+e)(1 - \bar{\xi}_s(x)) \right] ds.$$

By (4.4),

$$|\bar{\xi}_t| \leq |\bar{\xi}_0| + \int_0^t b|\bar{\xi}_s| \, ds + \bar{M}_t,$$

and as we have already noted that $|\bar{\xi}_t|$ has a finite mean, (4.2) follows by taking means in the above and using Gronwall's lemma.

Using some stochastic calculus in (4.4), we get (with $[M]_t$ the square variation function of $M_t$)

$$(4.6) \qquad |\bar{\xi}_t|^2 = |\bar{\xi}_0|^2 + \int_0^t 2|\bar{\xi}_s| b \sum_{x,e} \bar{p}(e) \bar{\xi}_s(x+e)(1 - \bar{\xi}_s(x)) \, ds$$
$$+ \int_0^t 2|\bar{\xi}_s| \, d\bar{M}_s + [\bar{M}]_t.$$

Proposition 2.1, the fact that $|\bar{\xi}_t|$ can be bounded by a pure birth process and (4.5) imply that the stochastic integral in the above is a martingale, as is $[\bar{M}]_t - \langle \bar{M} \rangle_t$, consequently,

$$(4.7) \qquad E(|\bar{\xi}_t|^2) \leq |\bar{\xi}_0|^2 + 2b \int_0^t E(|\bar{\xi}_s|^2) \, ds + E(\langle \bar{M} \rangle_t)$$
$$\leq |\bar{\xi}_0|^2 + 2b \int_0^t E(|\bar{\xi}_s|^2) \, ds + \int_0^t (2v+b) E(|\bar{\xi}_s|) \, ds.$$

From this, (4.2) and the previously noted fact that $E(|\xi_t|^2)$ is bounded on compact time intervals, (4.3) is easy to derive. Finally, the fact that $|\bar{\xi}_t|$ and $|\bar{\xi}_t|^2$ are submartingales is clear from (4.4) and (4.6). $\square$

Our task now is to define a biased voter model $\bar{\xi}_t^N$ taking values in $\{0,1\}^{\mathbf{S_N}}$ which dominates the voter model perturbation $\xi_t^N$. To do this, we must determine the appropriate kernels and rates $v = v_N$ and $b = b_N$, which we do by considering the maximum and minimum values of $c_N(x,\xi)$ given by (1.14), (1.15) and (1.17). We assume that $N \geq k_\delta$ [recall (P4)] in what follows.

For $\xi_t^N$, at site $x$ in configuration $\xi$ with $\xi(x) = 1$, the flip rate from 1 to 0 is

$$
\begin{aligned}
c_N(x,\xi) &= N \sum_{y \in \mathbf{S_N}} p_N(y-x)(1-\xi(y)) + \sum_{A \in P_F} \delta_N(A)\chi_N(A,x,\xi) \\
&\geq (N - k_\delta) f_0^N(x,\xi),
\end{aligned}
\tag{4.8}
$$

where we have made use of assumption (P4).

Similarly, at site $x$ in configuration $\xi$ with $\xi(x) = 0$, the flip rate from 0 to 1 is

$$
\begin{aligned}
c_N(x,\xi) &= N \sum_{y \in \mathbf{S_N}} p_N(y-x)\xi(y) + \sum_{A \in P_F} \beta_N(A)\chi_N(A,x,\xi) \\
&\leq N f_1^N(x,\xi) + \sum_{A \in P_F} \beta_N^+(A)\chi_N(A,x,\xi) \\
&\leq N f_1^N(x,\xi) + \sum_{A \in P_F} \frac{\beta_N^+(A)}{|A|} \sum_{a \in A} \xi(x + a/\ell_N),
\end{aligned}
\tag{4.9}
$$

where we have used (2.10). To simplify this last expression, we define a probability kernel $\hat{p}_N$ on $\mathbf{S_N}$ by setting $c_\beta^N = \sum_{A \in P_F} \beta_N^+(A)$ and

$$\hat{p}_N(a) = \frac{1}{c_\beta^N} \sum_{A : a \in A/\ell_N} \frac{\beta_N^+(A)}{|A|}.$$

(If $c_\beta^N = 0$, the construction simplifies considerably and the necessary modifications will be obvious.) Note that $\hat{p}_N(0) = 0$. Now if

$$\hat{f}_i^N(x,\xi) = \sum_y \hat{p}_N(y-x)\mathbb{1}\{\xi(y) = i\},$$

inequality (4.9) can be rewritten as

$$c_N(x,\xi) \leq N f_1^N(x,\xi) + c_\beta^N \hat{f}_1^N(x,\xi). \tag{4.10}$$

Recall by (P1), $c_\beta = \sup_N c_\beta^N < \infty$, and we use this constant to define another probability kernel $\bar{p}_N$ on $\mathbf{S_N}$ by

$$\bar{p}_N(a) = \frac{k_\delta p_N(a) + c_\beta \hat{p}(a)}{k_\delta + c_\beta}.$$

It follows then, with $\bar{f}_i^N(x,\xi) = \sum_y \bar{p}_N(y-x)\mathbb{1}\{\xi(y) = i\}$, that (recall $\bar{c} = k_\delta + c_\beta$)

$$(4.11) \qquad Nf_1^N(x,\xi) + c_\beta^N \hat{f}_1^N(\xi) \leq (N - k_\delta)f_1^N(x,\xi) + \bar{c}\bar{f}_1^N(\xi).$$

We now let $\bar{\xi}_t^N$ be the biased voter model with rate function

$$(4.12) \qquad \bar{c}_N(x,\xi) = \begin{cases} (N - k_\delta)f_1^N(x,\xi) + \bar{c}\bar{f}_1^N(\xi), & \text{if } \xi(x) = 0, \\ (N - k_\delta)f_0^N(x,\xi), & \text{if } \xi(x) = 1. \end{cases}$$

From (4.8), (4.10) and (4.11), we see that if $\xi \leq \bar{\xi}$,

$$(4.13) \qquad \begin{aligned} c_N(x,\xi) &\leq \bar{c}_N(x,\bar{\xi}) && \text{if } \bar{\xi}(x) = 0, \\ c_N(x,\xi) &\geq \bar{c}_N(x,\bar{\xi}) && \text{if } \xi(x) = 1. \end{aligned}$$

On account of this (see Theorem III.1.5 of [11]), we may construct versions of $\xi_t^N$ and $\bar{\xi}_t^N$ on a common probability space such that if $\xi_0^N = \bar{\xi}_0^N$, then with probability one,

$$(4.14) \qquad \xi_t^N \leq \bar{\xi}_t^N \qquad \text{for all } t \geq 0.$$

In Section 5 we will also need a voter model dominated by $\bar{\xi}_t^N$. Let $\hat{\xi}_t^N$ be the process with the same flip rates specified in (4.12), except with $\bar{c} = 0$. Then $\hat{\xi}_t^N$ is a voter model, and if $\hat{\xi}_0^N(x) \leq \bar{\xi}_0^N(x)$ for all $x$, then, as above, we can define $\hat{\xi}_t^N$ and $\bar{\xi}_t^N$ on a common probability space so that with probability one,

$$(4.15) \qquad \hat{\xi}_t^N \leq \bar{\xi}_t^N \qquad \text{for all } t \geq 0.$$

We also note that $|\hat{\xi}_t^N|$ is a *martingale* [e.g., by (4.4) with $b = 0$], so

$$(4.16) \qquad E(|\hat{\xi}_t^N|) = |\hat{\xi}_0^N| \qquad \text{for all } t \geq 0.$$

We record now some consequences of Lemma 4.1, including the proof of Proposition 3.3. We assume that $\bar{X}_t^N$ and $\hat{X}_t^N$ are as above, with $\bar{\xi}_0^N = \hat{\xi}_0^N = \xi_0^N$. Let $\bar{X}_t^N(\phi) = (1/N)\sum_x \phi(x)\bar{\xi}_t^N(x)$ and $\hat{X}_t^N(\phi) = (1/N)\sum_x \phi(x)\hat{\xi}_t^N(x)$. By Lemma 4.1,

$$(4.17) \qquad E(\bar{X}_t^N(\mathbf{1})) \leq e^{\bar{c}t}\bar{X}_0^N(\mathbf{1}).$$

Also by Lemma 4.1,

$$(4.18) \quad E\bar{X}_t^N(\mathbf{1})^2 \leq e^{2\bar{c}t}\left(\bar{X}_0^N(\mathbf{1})^2 + \frac{\bar{c} + 2(N-k_\delta)}{N\bar{c}}(1 - e^{-\bar{c}t})\bar{X}_0^N(\mathbf{1})\right).$$

Since $\bar{X}_t^N(\mathbf{1})^2$ is a submartingale by Lemma 4.1, it follows that for $T > 0$ and $K > 0$, there exists a constant $C(T,K) \geq 1$ such that

$$(4.19) \qquad \sup_{\bar{X}_0^N(\mathbf{1}) \leq K} E\left(\sup_{t \leq T} \bar{X}_t^N(\mathbf{1})^2\right) \leq C(T,K).$$

PROOF OF PROPOSITION 3.3. This is now immediate from the above inequality, since the coupling $\xi_t^N \leq \bar{\xi}_t^N$ implies that $X_t^N(\mathbf{1}) \leq \bar{X}_t^N(\mathbf{1})$. □

Note that by (4.17) and the fact that $\bar{X}_t^N(\mathbf{1})$ is a submartingale,

$$(4.20) \qquad 0 \leq E(\bar{X}_t^N(\mathbf{1})) - \bar{X}_0^N(\mathbf{1}) \leq (e^{\bar{c}t} - 1)\bar{X}_0^N(\mathbf{1}).$$

To get similar bounds on the difference $X_t^N(\mathbf{1}) - X_0^N(\mathbf{1})$, use Proposition 2.3 and Lemma 3.5 to see that $X_t^N(\mathbf{1}) - X_0^N(\mathbf{1}) = \int_0^t d_s^N(\mathbf{1})\,ds + M_t^N(\mathbf{1})$, where $E(M_t^N(\mathbf{1})) = 0$, and there is a constant $C$ such that

$$|d_s^N(\mathbf{1})| \leq CX_s^N(\mathbf{1}) \leq C\bar{X}_s^N(\mathbf{1})$$

for $s \leq T$. It follows therefore from (4.17) that

$$(4.21) \qquad |E(X_t^N(\mathbf{1}) - X_0^N(\mathbf{1}))| \leq C\frac{e^{\bar{c}t} - 1}{\bar{c}}X_0^N(\mathbf{1}).$$

**5. The key lemma.** For bounded functions $\phi$ on $\mathbf{S_N}$ and nonempty $A \in P_F$, define

$$(5.1) \quad \begin{aligned} &\eta_N(X_0^N, A, \phi, s) \\ &= \left|E_{X_0^N}\left(\frac{1}{N}\sum_x \phi(x)\chi_N(A, x, \xi_s^N) - P(\tau^N(A) \leq s)X_s^N(\phi)\right)\right| \end{aligned}$$

and

$$(5.2) \qquad \eta_{N,J}(A, \phi, s) = \sup_{X_0^N(\mathbf{1}) \leq J} \eta_N(X_0^N, A, \phi, s).$$

The proof of Proposition 3.4 is based on the following lemma. We assume the hypotheses of Theorem 1.3 are in force.

LEMMA 5.1. *There is a finite constant $C$ and a positive sequence $\varepsilon_N \to 0$ as $N \to \infty$ such that for any $J, K \geq 1$, $\phi\colon\mathbf{S_N} \to \mathbb{R}$ such that $\|\phi\|_{\mathrm{Lip}} \leq K$, nonempty finite $A \subset \mathbb{Z}^d$ and $\bar{a} \in A$, and $s > 0$,*

$$(5.3) \quad \begin{aligned} &\eta_N(X_0^N, A, \phi, s) \\ &\leq CK\left[(e^{\bar{c}s} - 1)|A| + \left(P(\tau^N(A) \leq s) \wedge \left(\frac{|\bar{a}|}{\ell_N} + E|B_s^{N,0}|\right)\right)\right]X_0^N(\mathbf{1}) \\ &\quad + CK|A|NP(B_s^{N,0} = 0)(X_0^N(\mathbf{1}))^2 \end{aligned}$$

*and*

$$(5.4) \qquad \eta_{N,J}(A, \phi, \varepsilon_N^*) \leq CKJ^2\left(\varepsilon_N|A| + \sigma_N(A) \wedge \left(\frac{|\bar{a}|}{\ell_N} + \varepsilon_N\right)\right).$$

PROOF. Let $J, K, \phi$ and $A$ be as above. Let $\bar{\xi}_t^N$ be the biased voter model and let $\hat{\xi}_t^N$ be the voter model from the previous section, with $\xi_0^N = \bar{\xi}_0^N = \hat{\xi}_0^N$, coupled so that $\xi_t^N \leq \bar{\xi}_t^N$ and $\hat{\xi}_t^N \leq \bar{\xi}_t^N$. By the triangle inequality, $\eta_N(X_0^N, A, \phi, s)$ is bounded above by the sum of the following four "error" terms:

$$(5.5) \quad \eta_1^N(s) = \left| E\left( \frac{1}{N} \sum_x \phi(x)[\chi_N(A, x, \xi_s^N) - \chi_N(A, x, \bar{\xi}_s^N)] \right) \right|,$$

$$(5.6) \quad \eta_2^N(s) = \left| E\left( \frac{1}{N} \sum_x \phi(x)[\chi_N(A, x, \bar{\xi}_s^N) - \chi_N(A, x, \hat{\xi}_s^N)] \right) \right|,$$

$$(5.7) \quad \eta_3^N(s) = \left| E\left( \left[ \frac{1}{N} \sum_x \phi(x)\chi_N(A, x, \hat{\xi}_s^N) \right] - P(\tau^N(A) \leq s)\hat{X}_0^N(\phi) \right) \right|,$$

$$(5.8) \quad \eta_4^N(s) = P(\tau^N(A) \leq s)|E(\hat{X}_0^N(\phi) - X_s^N(\phi))|$$

(recall $\hat{X}_0^N = X_0^N$).

The strategy behind this decomposition is as follows. We want to argue that for small $s$, the perturbed voter model $\xi_s^N$ is close in some sense to the voter model $\hat{\xi}_s^N$, and then compute with $\hat{\xi}_s^N$ using voter model *duality*. However, we cannot directly compare $\xi_s^N$ with $\hat{\xi}_s^N$, but must instead argue that both $\xi_s^N$ and $\hat{\xi}_s^N$ are close to $\bar{\xi}_s^N$. These two comparisons can be made because of the couplings and the inequality $|\prod_{i=1}^n z_i - \prod_{i=1}^n w_i| \leq \sum_{i=1}^n |z_i - w_i|$ for numbers $z_i, w_i$ bounded in absolute value by 1.

In preparation for estimating the $\eta_i^N(s)$, by the previous inequality,

$$|\chi_N(A, x, \xi_s^N) - \chi_N(A, x, \bar{\xi}_s^N)| \leq \sum_{a \in A} |\bar{\xi}_s^N(x + a/\ell_N) - \xi_s^N(x + a/\ell_N)|$$

$$= \sum_{a \in A} (\bar{\xi}_s^N(x + a/\ell_N) - \xi_s^N(x + a/\ell_N)),$$

the last step following from the coupling $\xi_s^N \leq \bar{\xi}_s^N$. Thus,

$$\frac{1}{N} \sum_{x \in \mathbf{S_N}} |\chi_N(A, x, \xi_s^N) - \chi_N(A, x, \bar{\xi}_s^N)| \leq |A|(\bar{X}_s^N(\mathbf{1}) - X_s^N(\mathbf{1})).$$

A similar argument shows that

$$\frac{1}{N} \sum_{x \in \mathbf{S_N}} |\chi_N(A, x, \bar{\xi}_s^N) - \chi_N(A, x, \hat{\xi}_s^N)| \leq |A|(\bar{X}_s^N(\mathbf{1}) - \hat{X}_s^N(\mathbf{1})).$$

Consider the first error term $\eta_1^N(s)$. By the above,

$$\eta_1^N(s) \leq \frac{1}{N} \sum_{x \in \mathbf{S_N}} |\phi(x)|E|\chi_N(A, x, \bar{\xi}_s^N) - \chi_N(A, x, \xi_s^N)|$$

$$\leq \|\phi\|_\infty |A| E(\bar{X}_s^N(\mathbf{1}) - X_s^N(\mathbf{1}))$$
$$\leq K|A|(E(\bar{X}_s^N(\mathbf{1}) - \bar{X}_0^N(\mathbf{1})) + |E(X_0^N(\mathbf{1}) - X_s^N(\mathbf{1}))|)$$

[recall $\bar{X}_0^N(\mathbf{1}) = X_0^N(\mathbf{1})$]. By (4.20) and (4.21), this implies there is a constant $C$ such that

(5.9) $$\eta_1^N(s) \leq CK(e^{\bar{c}s} - 1)|A|X_0^N(\mathbf{1}).$$

For $\eta_2^N(s)$, using $E(\hat{X}_s^N(\mathbf{1})) = \hat{X}_0^N(\mathbf{1}) = X_0^N(\mathbf{1})$ [see (4.16)] and arguing as above, we get

$$\eta_2^N(s) \leq \|\phi\|_\infty |A| E(\bar{X}_s^N(\mathbf{1}) - \hat{X}_s^N(\mathbf{1})) \leq K|A|E(\bar{X}_s^N(\mathbf{1}) - X_0^N(\mathbf{1})).$$

Now apply (4.20) to see there is a constant $C$ such that

(5.10) $$\eta_2^N(s) \leq CK|A|(e^{\bar{c}s} - 1)X_0^N(\mathbf{1}).$$

Turning to $\eta_4^N(s)$, by adding and subtracting $\bar{X}_s^N(\phi)$ and then proceeding as above, there is a constant $C$ such that

(5.11) $$\eta_4^N(s) \leq CK(e^{\bar{c}s} - 1)X_0^N(\mathbf{1}).$$

We come now to the main term, $\eta_3^N(s)$. Here we will use the independent random walk system $\{B_t^{N,x}, x \in \mathbf{S_N}\}$ and the coalescing random walk system $\{\hat{B}_t^{N,x}, x \in \mathbf{S_N}\}$ introduced in Section 1. Recall that for $A \in P_F$,

$$\tau^N(A) = \inf\{t : |\{\hat{B}_t^{N,x}, x \in A/\ell_N\}| = 1\}.$$

For $y \in \mathbf{S_N}$, let $\tau_y^N(A) = \tau^N(y\ell_N + A)$. By translation invariance and symmetry, for any $y \in \mathbf{S_N}$ and finite $A \subset \mathbb{Z}^d$,

(5.12) $$P(\tau_y^N(A) \leq s) = P(\tau_0^N(A) \leq s) = P(\tau_0^N(-A) \leq s) = P(\tau_y^N(-A) \leq s).$$

Also, we may assume here that our coalescing random walk system is constructed from the independent random walk system via some collision rule. In particular, for $a \neq a' \in \mathbb{Z}^d$, we may assume that

(5.13) $$\begin{aligned} &P(\hat{B}_s^{N,x+a/\ell_N} = y, \hat{B}_s^{N,x+a'/\ell_N} = z, \tau_x^N(\{a,a'\}) > s) \\ &= P(B_s^{N,x+a/\ell_N} = y, B_s^{N,x+a'/\ell_N} = z, \tau_x^N(\{a,a'\}) > s) \\ &\leq P(B_s^{N,x+a/\ell_N} = y)P(B_s^{N,x+a'/\ell_N} = z). \end{aligned}$$

Finally, we will make use of the well-known duality between the voter model and coalescing random walk (see Section 3 of [5], e.g.) in the form

(5.14) $$E(\chi_N(A, x, \hat{\xi}_s^N)) = P(\hat{B}_s^{N,x+a/\ell_N} \in \hat{\xi}_0^N \,\forall\, a \in A).$$

We will evaluate the right-hand side above by decomposing the event according to whether $\tau_x^N(A) \leq s$ or not.

To estimate $\eta_3^N(s)$, we define

$$(5.15) \quad \eta_{3,1}^N(s) = \left| \frac{1}{N} \sum_{x \in \mathbf{S_N}} \phi(x) P(\hat{B}_s^{N,x+a/\ell_N} \in \hat{\xi}_0^N \,\forall a \in A, \tau_x^N(A) > s) \right|,$$

$$(5.16) \quad \eta_{3,2}^N(s) = \left| \frac{1}{N} \sum_{x \in \mathbf{S_N}} \phi(x) P(\hat{B}_s^{N,x+a/\ell_N} \in \hat{\xi}_0^N \,\forall a \in A, \tau_x^N(A) \le s) \right.$$

$$\left. - P(\tau^N(A) \le s) \hat{X}_0^N(\phi) \right|,$$

and observe that the duality equation (5.14) above implies that

$$\eta_3^N(s) \le \eta_{3,1}^N(s) + \eta_{3,2}^N(s).$$

We proceed now to estimate each of these terms.

For $\eta_{3,1}^N(s)$, fix any $\bar{a} \in A$. Since $\{\tau_x^N(A) > s)\} = \bigcup_{a \in A \setminus \{\bar{a}\}} \{\tau_x^N(a, \bar{a}) > s\}$, it follows from (5.13) and $P(B_s^{N,z} = w) \le P(B_s^{N,0} = 0)$ (e.g., see Lemma A.3 of [3]) that

$$\eta_{3,1}^N(s) \le \|\phi\|_\infty \frac{1}{N} \sum_{x \in \mathbf{S_N}} \sum_{a \in A \setminus \{\bar{a}\}} P(\hat{B}_s^{N,x+a/\ell_N} \in \hat{\xi}_0^N,$$

$$\hat{B}_s^{N,x+\bar{a}/\ell_N} \in \hat{\xi}_0^N, \tau_x^N(a,\bar{a}) > s)$$

$$\le \|\phi\|_\infty \sum_{a \in A \setminus \{\bar{a}\}} \frac{1}{N} \sum_{x,y,z \in \mathbf{S_N}} P(B_s^{N,x+a/\ell_N} = y)$$

$$\times P(B_s^{N,x+\bar{a}/\ell_N} = z) \hat{\xi}_0^N(y) \hat{\xi}_0^N(z)$$

$$\le \|\phi\|_\infty \sum_{a \in A \setminus \{\bar{a}\}} \frac{1}{(N)^2} \sum_{x,y,z \in \mathbf{S_N}} N P(B_s^{N,0} = 0)$$

$$\times P(B_s^{N,x+\bar{a}/\ell_N} = z) \hat{\xi}_0^N(y) \hat{\xi}_0^N(z).$$

By symmetry and time reversal, $P(\hat{B}_s^{N,x+\bar{a}/\ell_N} = z) = P(\hat{B}_s^{N,z} = x + \bar{a}/\ell_N)$. Thus, in the inequality above, if we carry out the summation first over $x$, and then over $y$ and $z$, we obtain the estimate

$$(5.17) \quad \eta_{3,1}^N(s) \le K(|A|-1) N P(B_s^{N,0} = 0)(\hat{X}_0^N(\mathbf{1}))^2.$$

For $\eta_{3,2}^N(s)$, we begin with a calculation that uses time reversal, symmetry and translation invariance. For any $\bar{a} \in A$,

$$P(\hat{B}_s^{N,x+\bar{a}/\ell_N} = y, \tau_x^N(A) \le s)$$

$$= P(\hat{B}_s^{N,0} = y - (x + \bar{a}/\ell_N), \hat{B}_s^{N,0} = \hat{B}_s^{N,(a-\bar{a})/\ell_N} \,\forall a \in A)$$

$$= P(\hat{B}_s^{N,0} = (x + \bar{a}/\ell_N) - y, \hat{B}_s^{N,0} = \hat{B}_s^{N,(\bar{a}-a)/\ell_N} \,\forall\, a \in A)$$
$$= P(\hat{B}_s^{N,y-\bar{a}/\ell_N} = x, \hat{B}_s^{N,y-\bar{a}/\ell_N} = \hat{B}_s^{N,y-a/\ell_N} \,\forall\, a \in A)$$
$$= P(\hat{B}_s^{N,y-\bar{a}/\ell_N} = x, \tau_y^N(-A) \leq s).$$

Using this equality, we have, for any fixed $\bar{a} \in A$,

$$\frac{1}{N} \sum_{x \in \mathbf{S_N}} \phi(x) P(\hat{B}_s^{N,x+a/\ell_N} \in \hat{\xi}_0^N \,\forall\, a \in A, \tau_x^N(A) \leq s)$$
$$= \frac{1}{N} \sum_{x,y \in \mathbf{S_N}} \phi(x) \hat{\xi}_0^N(y) P(\hat{B}_s^{N,x+\bar{a}/\ell_N} = y, \tau_x^N(A) \leq s)$$
$$= \frac{1}{N} \sum_{x,y \in \mathbf{S_N}} \phi(x) \hat{\xi}_0^N(y) P(\hat{B}_s^{N,y-\bar{a}/\ell_N} = x, \tau_y^N(-A) \leq s)$$
$$= \frac{1}{N} \sum_{y \in \mathbf{S_N}} \hat{\xi}_0^N(y) E(\phi(\hat{B}_s^{N,y-\bar{a}/\ell_N}); \tau_y^N(-A) \leq s).$$

Furthermore, since $P(\tau^N(A) \leq s) = P(\tau_y^N(-A) \leq s)$ for all $y \in \mathbf{S_N}$ [by (5.12)], adding and subtracting $\phi(y)$ in the sum above gives

$$\frac{1}{N} \sum_{x \in \mathbf{S_N}} \phi(x) P(\hat{B}_s^{N,x+a/\ell_N} \in \hat{\xi}_0^N \,\forall\, a \in A, \tau_x^N(A) \leq s)$$
$$= \frac{1}{N} \sum_{y \in \mathbf{S_N}} \hat{\xi}_0^N(y) E(\phi(\hat{B}_s^{N,y-\bar{a}/\ell_N}) - \phi(y); \tau_y^N(-A) \leq s)$$
$$+ P(\tau^N(A) \leq s) \hat{X}_0^N(\phi).$$

Therefore,

$$\eta_{3,2}^N(s) \leq \frac{1}{N} \sum_{y \in \mathbf{S_N}} \hat{\xi}_0^N(y) |E(\phi(\hat{B}_s^{N,y-\bar{a}/\ell_N}) - \phi(y); \tau_y^N(-A) \leq s)|.$$

Now, since $\|\phi\|_{\mathrm{Lip}} \leq K$,

$$|E(\phi(B_s^{N,y-\bar{a}/\ell_N}) - \phi(y); \tau_y^N(-A) \leq s)|$$
$$\leq (2KP(\tau^N(A) \leq s)) \wedge E\left(\left|\phi\left(y - \frac{\bar{a}}{\ell_N} + B_s^{N,0}\right) - \phi(y)\right|\right)$$
$$\leq 2K\left(P(\tau^N(A) \leq s) \wedge \left(\frac{|\bar{a}|}{\ell_N} + E(|B_s^{N,0}|)\right)\right).$$

Assembling these estimates, we obtain

$$\eta_{3,2}^N(s) \leq 2K \hat{X}_0^N(\mathbf{1})\left(P(\tau^N(A) \leq s) \wedge \left(\frac{|\bar{a}|}{\ell_N} + E(|B_s^{N,0}|)\right)\right).$$

It now follows from the estimates on $\eta_{3,1}^N(s)$ and $\eta_{3,2}^N(s)$ that

$$
\eta_3^N(s) \leq 2KX_0^N(\mathbf{1})\left[P(\tau^N(A) \leq s) \wedge \left(\frac{|\bar{a}|}{M_N\sqrt{N}} + E|B_s^{N,0}|\right)\right]
$$
(5.18)
$$
+ K|A|NP(B_s^{N,0} = 0)(X_0^N(\mathbf{1}))^2.
$$

Combining (5.9)–(5.11) and (5.18) completes the proof of (5.3). Setting $s = \varepsilon_N^*$ in (5.3) and using the kernel assumption (K1), we obtain (5.4), provided that $E(|B_{\varepsilon_N^*}^{N,0}|) \to 0$ as $N \to \infty$. But this follows easily from (H1), since $E(|B_{\varepsilon_N^*}^{N,0}|^2) = \varepsilon_N^* E(|W_N|^2)$. $\square$

**6. Proof of Proposition 3.4.** Let $T, K, \phi, A, J$ and $0 \leq t \leq T$ be as in the statement of Proposition 3.4. Define the hitting times

$$T_J^N = \inf\{s \geq 0 : X_s^N(\mathbf{1}) > J\}.$$

By Proposition 3.3,

(6.1)
$$\sup_N P(T_J^N \leq t) \leq C_3(K,T)J^{-2}.$$

Let $\varepsilon_N^* > 0$ be as in (K1)–(K3). Also, define

$$\Delta^N(A, \phi_s, \xi_s^N) = \frac{1}{N}\sum_x \phi_s(x)\chi_N(A, x, \xi_s^N) - P(\tau_N(A) \leq \varepsilon_N^*)X_s^N(\phi_s).$$

*Step* 1. We claim that for $t \leq T$,

$$
E\left(\left(\int_{(T_J^N + \varepsilon_N^*)\wedge t}^t \Delta^N(A, \phi_s, \xi_s^N)\,ds\right)^2\right)
$$
(6.2)
$$
\leq 4K^2 T C_3(K,T) J^{-2} \int_0^t E(X_s^N(\mathbf{1}))^2\,ds.
$$

This inequality is easily derived. For any $\bar{a} \in A$, $\chi_N(A, x, \xi_s^N) \leq \xi_s^N(x + \bar{a}/\ell_N)$, and, hence,

$$
|\Delta^N(A, \phi_s, \xi_s^N)| \leq \frac{1}{N}\sum_x |\phi_s(x)|\left(\xi_s^N\left(x + \frac{\bar{a}}{\ell_N}\right) + \xi_s^N(x)\right)
$$
(6.3)
$$
\leq 2\|\phi\|_\infty X_s^N(\mathbf{1}).
$$

With this inequality, Cauchy–Schwarz implies

$$
\left(\int_0^t \mathbb{1}\{s > T_J^N + \varepsilon_N^*\}\Delta^N(A, \phi_s, \xi_s^N)\,ds\right)^2
$$
$$
\leq tP(T_J^N \leq t)4\|\phi\|_\infty^2 \int_0^t (X_s^N(\mathbf{1}))^2\,ds
$$

and the claim follows from (6.1).

*Step* 2. Because $\mathbb{1}\{T_J^N < s_1 < T_J^N + \varepsilon_N^*\}\mathbb{1}\{s_1 + \varepsilon_N^* < s_2 < T_J^N + \varepsilon_N^*\} = 0$,

(6.4) $\quad E\bigg(\bigg(\int_0^{(T_J^N+\varepsilon_N^*)\wedge t} \Delta^N(A,\phi_s,\xi_s^N)\,ds\bigg)^2\bigg) = I_1(N,J,t) + I_2(N,J,t),$

where

(6.5)
$$I_1(N,J,t) = 2\int_0^t E\bigg[\mathbb{1}_{\{s_1 \leq T_J^N + \varepsilon_N^*\}}\Delta^N(A,\phi_{s_1},\xi_{s_1}^N)$$
$$\times \int_{s_1}^{(s_1+\varepsilon_N^*)\wedge t}\mathbb{1}_{\{s_2 \leq T_J^N+\varepsilon_N^*\}}\Delta^N(A,\phi_{s_2},\xi_{s_2}^N)\bigg)ds_2\bigg]ds_1$$

and

(6.6)
$$I_2(N,J,t) = 2\int_0^t E\bigg[\mathbb{1}_{\{s_1 \leq T_J^N\}}\Delta^N(A,\phi_{s_1},\xi_{s_1}^N)$$
$$\times \int_{(s_1+\varepsilon_N^*)\wedge t}^t \mathbb{1}_{\{s_2 \leq T_J^N + \varepsilon_N^*\}}\Delta^N(A,\phi_{s_2},\xi_{s_2}^N)\,ds_2\bigg]ds_1.$$

By (6.3), (4.17) and the Markov property,

$$|I_1(N,J,t)| \leq 8\|\phi\|_\infty^2 E\bigg(\int_0^t X_{s_1}^N(\mathbf{1})\int_{s_1}^{s_1+\varepsilon_N^*} X_{s_2}^N(\mathbf{1})\,ds_2\,ds_1\bigg)$$
$$\leq 8\|\phi\|_\infty^2 E\bigg(\int_0^t X_{s_1}^N(\mathbf{1})\int_{s_1}^{s_1+\varepsilon_N^*} E_{X_{s_1}^N}(X_{s_2-s_1}^N(\mathbf{1}))\,ds_2\,ds_1\bigg)$$
$$\leq 8\|\phi\|_\infty^2 E\bigg(\int_0^t X_{s_1}^N(\mathbf{1})\varepsilon_N^* e^{\bar{c}\varepsilon_N^*} X_{s_1}^N(\mathbf{1})\,ds_1\bigg)$$
$$= 8K^2\varepsilon_N^* e^{\bar{c}\varepsilon_N^*}\int_0^t E(X_s^N(\mathbf{1}))^2\,ds.$$

Now consider $I_2(N,J,t)$. Let $0 \leq s_1 < s_2 < t$ satisfy $s_1 + \varepsilon_N^* < s_2 < T_J^N + \varepsilon_N^*$, in which case $X_{s_2-\varepsilon_N^*}^N(\mathbf{1}) \leq J$. Then

$$|E(\mathbb{1}\{s_1 < T_J^N\}\mathbb{1}\{s_2 < T_J^N + \varepsilon_N^*\}\Delta_N(A,\phi_{s_1},\xi_{s_1}^N)\Delta_N(A,\phi_{s_2},\xi_{s_2}^N))|$$
$$\leq E(\mathbb{1}\{s_1 < T_J^N\}\mathbb{1}\{s_2 < T_J^N + \varepsilon_N^*\}|\Delta_N(A,\phi_{s_1}\xi_{s_1}^N)|$$
$$\times |E_{X_{s_2-\varepsilon_N^*}^N}(\Delta_N(A,\phi_{s_2},\xi_{\varepsilon_N^*}^N)|))$$
$$\leq E(\mathbb{1}\{s_1 < T_J^N\}\mathbb{1}\{s_2 < T_J^N + \varepsilon_N^*\}|\Delta_N(A,\phi_{s_1},\xi_{s_1}^N)|\eta_{N,J}(A,\phi_{s_2},\varepsilon_N^*))$$
$$\leq \eta_{N,J}(A,\phi_{s_2},\varepsilon_N^*)2\|\phi\|_\infty E(X_{s_1}^N(\mathbf{1})),$$

the last by (6.3). By these estimates we have

(6.7) $\quad I_2(N,J,t) \leq 2\int_0^t \eta_{N,J}(A,\phi_s,\varepsilon_N^*)\,ds\, 2K\int_0^t E(X_s^N(\mathbf{1}))\,ds.$

Now for the proof of (3.2). By the above bounds, and Proposition 3.3 and Lemma 5.1, if $\bar{a}$, $\varepsilon_N$ and $J$ are as in Lemma 5.1, then for $t \leq T$,

$$\mathcal{E}_N(A, \phi, K, t)$$
$$\leq C(K,T) \bigg[ (J^{-2} + \varepsilon_N^* e^{\bar{c}\varepsilon_N^*}) \int_0^T E(X_s^N(\mathbf{1})^2) \, ds$$
$$+ \int_0^T \eta_{N,J}(A, \phi_s, \varepsilon_N^*) \, ds \int_0^T E(X_s^N(\mathbf{1})) \, ds \bigg]$$
$$\leq C(K,T) \bigg[ J^{-2} + \varepsilon_N^* e^{\bar{c}\varepsilon_N^*} + J^2 \bigg( \varepsilon_N |A| + \sigma_N(A) \wedge \bigg( \frac{\bar{a}}{\ell_N} + \varepsilon_N \bigg) \bigg) \bigg],$$

and we are done.

**Acknowledgment.** E. A. Perkins thanks Alison Etheridge for some interesting conversations on spatial stochastic Lotka–Volterra equations in mathematical ecology, and for telling us of her parallel work with Blath and Meredith [1] on co-existence of types using a different sde-type stochastic model.

DEPARTMENT OF MATHEMATICS
SYRACUSE UNIVERSITY
SYRACUSE, NEW YORK 13244
USA
E-MAIL: jtcox@syr.edu

DEPARTMENT OF MATHEMATICS
UNIVERSITY OF BRITISH COLUMBIA
VANCOUVER, BRITISH COLUMBIA
CANADA
E-MAIL: perkins@math.ubc.ca